\documentclass[final,3p,times]{elsarticle}
\usepackage{amssymb}
\usepackage{amsmath} 
\usepackage{color}
\usepackage{algorithm}
\usepackage{algpseudocode}
\usepackage{float}

\algnewcommand{\algorithmicforeach}{\textbccf{for each}}
\algdef{SE}[FOR]{ForEach}{EndForEach}[1]
  {\algorithmicforeach\ #1\ \algorithmicdo}% \ForEach{#1}
  {\algorithmicend\ \algorithmicforeach}% \EndForEach

\begin{document}

\begin{frontmatter}



\title{Data-driven multifidelity topology design with multi-channel variational auto-encoder for concurrent optimization of multiple design variable fields}
%\title{Concurrent optimization of high-fidelity modeling parameter and material distribution for data-driven multifidelity topology design}


\author[a]{Hiroki Kawabe\corref{cor1}}
\ead{kawabe@syd.mech.eng.osaka-u.ac.jp}
\cortext[cor1]{Corresponding author.}
\author[a]{Kentaro Yaji}
\author[b]{Yuichiro Aoki}

\address[a]{Department of Mechanical Engineering, Graduate School of Engineering, Osaka University, 2-1, Yamadaoka, Suita, Osaka 565-0871, Japan}
\address[b]{Aviation Technology Directorate, Japan Aerospace Exploration Agency, 6-13-1, Mitaka, Tokyo, 181-0015, Japan}

\begin{abstract}
% トポロジー最適化の概要
Topology optimization can generate high-performance structures with a high degree of freedom.
% トポロジー最適化の課題
Regardless, it generally confronts entrapment in undesirable local optima especially in problems characterized by strong non-linearity.
%However, conventional data-driven topology optimization methods have limitations in exploring modeling parameters like thickness to extrude the material distribution, which are crucial for identifying high-performance structures.
This study aims to establish a gradient-free topology optimization framework that facilitates more global solution searches to avoid the entrapment.
The framework utilizes a data-driven multifidelity topology design (MFTD), where solution candidates initially generated by solving low-fidelity (LF) optimization problems are iteratively updated by a variational auto-encoder (VAE) and high-fidelity (HF) evaluation.
%In this study, we propose a framework that enables the concurrent optimization of multiple design variable fields such as material and thickness distributions, utilizing data-driven multifielity topology design (MFTD) whose basic idea is to divide the original optimization problem into two components: low-fidelity (LF) optimization and high-fidelity (HF) evaluation.
A key procedure of the solution update is to construct HF models by extruding material distributions obtained by the VAE to a specified thickness called HF modeling parameter, which is conventionally constant across all solution candidates.
This constant assignment leads to no exploration of the HF modeling parameter space, which necessitates extensive parametric studies outside the optimization loop.
To enable a more comprehensive optimization in a single run, we propose a multi-channel image data architecture that stores material distributions in the first channel and HF modeling parameters in the second or subsequent channels.
%The framework is based on a multi-channel image data architecture, where the material distributions are stored in the first channel, and the HF modeling parameters are stored in the second or subsequent channels, utilizing multi-channel variational auto-encoder (MC-VAE).
This significant shift enables a thorough exploration of the HF modeling parameter space with no necessity of parametric studies afterwards, by simultaneously optimizing both material distributions and HF modeling parameters.
We apply the framework to a maximum stress minimization problem, where the LF optimization problem is formulated with approximation techniques, whereas the HF evaluation is conducted by accurately analyzing the stress field, bypassing any approximation techniques.
We first validate that the framework can successfully identify high-performance solutions superior to the reference solutions by effectively exploring both material distributions and HF modeling parameters in a fundamental stiffness maximization.
Then we demonstrate the framework can identify promising solutions for the original maximum stress minimization problems.
%This demonstration exhibits a distinct advantage in globally searching for optima across material distribution and modeling parameter spaces compared to the conventional data-driven MFTD, inticating the significant reduction of user reliance on finding optimal solutions by enabling a more comprehensive optimization in a single run, unlike conventional methods that require extensive parametric studies for various modeling parameters.

\end{abstract}

% =============後から追加=============
% =============後から追加=============


\end{frontmatter}

%\linenumbers

%sssssssssssssssssssssssssssssssssssssssssssssssssssssssssssssssssssssssssssssss
\section{Introduction}
\label{sec:intr}
% トポロジー最適化の概念的な紹介
% 構造設計におけるトポロジー最適化：最大応力最小化問題への導入
% 最大応力最小化問題の課題を列挙：強非線形性から勾配型の限界
% D2MFTDの導入：大局的に探索できる非勾配型トポロジー最適化、近似解でなく厳密解を評価できる
% D2MFTDの抱える課題
% 本研究における提案手法
% 本論文の構成
%sssssssssssssssssssssssssssssssssssssssssssssssssssssssssssssssssssssssssssssss

% トポロジー最適化の概念紹介
In structural design, there is an increasing demand for designing lightweight structures that are optimized for performances such as stiffness and strength.
One of the most successful tools to meet the demand is topology optimization, the concept of which is based on a formulation that shifts the focus from optimizing structural configurations to optimizing material distribution within the design domain, as introduced by Bends{\o}e and Kikuchi \cite{bendsoe1988generating}.
This elegantly simple yet effective formulation renders topology optimization sufficiently versatile for addressing a wide range of design problems, such as structural design and thermal-fluid design \cite{sigmund2013topology, deaton2014survey}.

% 構造設計におけるトポロジー最適化：剛性最大化問題への適用が進んでいるが、より重要性の高い応力集中や座屈といった非線形性の強い問題に対しては適用が進んでいない。
In practical use, topology optimization has been predominantly utilized for stiffness maximization problems that are one of the most fundamental optimization problem in structural design \cite{zhu2016topology}.
Despite its widespread application in such contexts, topology optimization has been less frequently applied to problems characterized by strong non-linearity, such as maximum stress minimization and buckling load maximization, compared to its use in stiffness maximization.
Notably, maximum stress minimization is crucial for designing structures within industrial applications, where failure criteria are primarily stress-based.
Consequently, numerous researchers are actively tackling the challenges associated with maximum stress minimization problems \cite{le2010stress, holmberg2013stress, duysinx1998topology}.

% 最大応力最小化問題の課題を列挙：非線形性が強く局所解に陥りやすい・疑似問題に帰着しており正確な応力の評価ができていない
Maximum stress minimization has a notable challenge called ``singularity phenomenon.''
This challenge arises when material removal leads to rapid increases in stress, causing numerical instabilities.
To mitigate this, traditional strategies have adopted stress-relaxation techniques.
Moreover, because discrete maximum stress values are impractical for optimization, these methods often resort to approximation techniques, substituting them with more continuous functions.
These approximations, however, introduce strong non-linearities, and in the prevalent approach of gradient-based topology optimization, there is a tendency to converge to local minima under these conditions.
Consequently, traditional approaches, through their reliance on relaxation and approximation techniques, can be viewed as addressing a modified version of the original optimization problem, rather than directly solving itself.

% 局所解を避けて大域的に探索できる非勾配型トポロジー最適化としてEAの使用を挙げ、その問題点を指摘
To avoid entrapment in local optima and promote a more global search, gradient-free topology optimization methods, such as Evolutionary Algorithms (EAs) \cite{coello2007evolutionary}, could offer a viable approach for addressing the challenges posed by strong non-linearities, which motivated a number of researchers to explore the application of EAs to topology optimization \cite{hu2010topology, madeira2010ga, munk2015topology}.
Although EAs enable a broader exploration of the solution space, the scalability of design variables in EAs is constrained to the order of $10^2$ to $10^3$, which is markedly less than the capacity of conventional gradient-based topology optimization methods, typically extending to the order of $10^5$ to $10^6$ \cite{sigmund2011usefulness}.
This constraint arises because two primary procedures of EAs, i.e., crossover and mutation, conventionally handle design variables as individual entities with no consideration of their spatial relationships.
Handling such a large number of design variables individually leads to the \textit{curse of dimensionality}, where the EA operations require an exponentially increasing computational cost as the number of design variables increases.

% Deep generative modelによるdata-driven design
To address the scalability issue of design variables in topology optimization, deep generative models \cite{foster2019generative} can be employed to reduce the dimension of the design space with a degree of freedom maintained.
In the realm of engineering design, data-driven designs based on those deep generative models has been actively investigated recently \cite{regenwetter2022deep, lee2024data}.
Guo et al. pioneered a data-driven topology optimization method that employs a Variational Auto-Encoder (VAE) \cite{kingma2013autoencoding} to encode and decode material distributions in the latent space, demonstrating the effectiveness of deep generative models in representing the topologically optimized designs \cite{guo2018indirect}.
Oh et al. introduced a Generative Adversarial Network (GAN) \cite{goodfellow2020generative} to generate diverse material distributions, which successfully identified high-performance structures \cite{oh2019deep}.
Wang et al. utilized a Latent Variable Gaussian Process (LVGP) to embed a mixed-variable inputs consisting of both qualitative microstructure concepts and quantitative microstructure design variables into a continuous and differentiable design space, enabling the multiscale topology optimization \cite{wang2020data}.

% EA-likeな操作を用いたD2MFTDの導入
To bypass the \textit{curse of dimensionality} in EA-based topology optimization, a series of researches have focused on a approach to replace the conventional crossover and mutation procedures with EA-like ones that inherit their essential concept but are specifically tuned for topology optimization requiring numerous design variables.
Yamasaki et al. proposed a data-driven topology design based on a crossover-like operation using a VAE that can compress high-dimensional input data into a low-dimensional latent space and then reconstruct high-dimensional output from this latent space.
They demonstrated that deep generative model-based crossover is effective to efficiently generate offspring material distributions that inherit the topological features from the parent ones \cite{yamasaki2021data}.
Yaji et al. introduced a mutation-like operation that deduces a mutant material distribution by solving topology optimization additionally constrained with referential material distribution in the dataset under optimization \cite{yaji2022data}.
They confirmed that the mutation-like operation is effective to prevent premature convergence in the iterative process through several numerical examples.

% MFTDの導入・利点
Moreover, Yaji et al. integrated this EA-based topology optimization method with Multi-Fidelity Topology Design (MFTD), based on the concept of multifidelity method known for its success in indirectly solving the original problem by addressing both low-fidelity (LF) and high-fidelity (HF) problems \cite{forrester2007multi}.
In MFTD, the original topology optimization problem is solved through two primary steps: LF optimization and HF evaluation.
The beauty of this multifidelity formulation is that a topology optimization problem, analytically or numerically difficult to directly solve, can be indirectly solved by replacing itself with a simplified optimization problem, then ensuring mechanical performance accurately evaluated through the optimization process.
The MFTD's effectiveness has been evidently confirmed by a vast array of mechanical applications, such as heat exchangers and redox flow batteries \cite{yaji2019framework, yaji2020multifidelity, kobayashi2019freeform}.

% 最大応力最小化問題に対するD2MFTDの適用例：Kato et al. (2023), Kii et al. (2024)
To tackle the strong non-linearity and the inaccurate stress evaluation by approximation techniques in maximum stress minimization problems, Kato et al. applied the data-driven MFTD to a maximum stress minimization problem \cite{kato2023tackling}.
They confirmed the accurate stress evaluation and more global solution search of the data-driven MFTD resulted in material distributions superior to those derived from the conventional gradient-based topology optimization.
Additionally, Kii et al. proposed a crossover-like operation called \textit{latent crossover} \cite{kii2024latent} that enhances the convergence of the data-driven MFTD by generating offspring material distributions with latent vectors sampled using the simplex crossover \cite{tsutsui1999multi} to inherit the parental characteristics more efficiently than the conventional random latent vectors uniformly sampled from the latent space.

% 最直近の文献の欠点を挙げる
% モデリングの方法について具体例を挙げて詳細に説明する
% 従来のD2MFTDでは、HF空間はパラメトリックに探索する、すなわち設計者への依存性が残っていた。
Despite the significant progress in data-driven MFTD for maximum stress minimization, the conventional data-driven MFTD has a limitation in exploring parameters for constructing the HF model.
In previous data-driven MFTD approaches, to construct the HF model, the thickness of the material distribution is typically assigned as a constant value across all solution candidates, with no exploration of the thickness value during the optimization process.
However, the thickness of the material distribution, we refer to as HF modeling parameters, has a significant impact on structural performance metrics, such as stiffness and stress, for the entire structure.
Therefore, exploring the HF parameter space in addition to the material distribution space is crucial for identifying high-performance structures.
To explore the HF parameter space, the conventional data-driven MFTD requires extensive parametric studies, which often highly depend on user intuition.

% 本研究では、HF空間の探索における依存性を排除し、その上で初期解に対して優れた解が得られることを確認する。
In this study, to expand the optimization search to the HF modeling space, we focus on an architecture of image data used in the data-driven MFTD.
The data-driven MFTD conventionally adopts a single-channel image to store material distributions, where HF modeling parameters are exclusively assigned as a constant value across each material distribution.
Instead of representing material distribution as a single-channel image, we introduce a multi-channel image where design variables for constructing the HF model are stored in the second or subsequent channels of multi-channel images, with the material distributions maintained in the first channel.
This significant shift enables a thorough exploration of the HF modeling parameter space with no necessity of parametric studies afterwards, by simultaneously optimizing both material distributions and HF modeling parameters.

% Multi-channelにしたことで必要になった変更について：mutationどこまで言及するか要検討
To transition the data architecture in the data-driven MFTD from a single-channel to multi-channel image, we primarily consider two aspects for adaptation: the crossover operation facilitated by VAE and the mutation operation achieved through solving constrained LF optimization.
Unlike the original data-driven MFTD, where the VAE's machine learning model architecture consisted simply of dense layers, this study incorporates a convolutional neural network (CNN) to efficiently learn the relationship between material distributions in the first channel and the HF parameters in subsequent channels \cite{gu2018recent}.
For the mutation operation, material distributions are generated in a manner similar to the traditional data-driven MFTD and then stored in the first channel of the multi-channel images.
After preparing the first channel, subsequent channels are designated to globally represent the HF parameter space.

% 2D to 3D mapping：論点がズレるのでDDMの次元拡張には触れない?
%Furthermore, concerning the DDM between the design domain and the unit plane structured mesh, traditional D2MFTD performed this mapping within a two-dimensional design domain. However, the authors suggest that this mapping technique could be extended to three-dimensional space. Given that structures requiring HF space exploration, such as thickness optimization, are often modeled as thin finite element (FE) models using three-dimensional shell elements, this study extends the mapping from two-dimensional to three-dimensional space.

% 実用的な構造設計の代表例として、補強外板構造を挙げつつ、HF空間の探索の重要性を説く
The primary goal of this study is to develop an enhanced data-driven MFTD capable of identifying high-performance structures within the HF modeling parameter space in addition to traditional material distribution space, especially for optimization problems characterized by significant non-linearity, such as maximum stress minimization.
Initially, we apply this approach to stiffness maximization for a simple reinforced skin structure, where the thickness of the reinforcement plays a crucial role in determining structural performance metrics, such as stiffness and stress.
Subsequently, we demonstrate the effectiveness of the proposed framework in addressing maximum stress minimization using the same structural model for the stiffness maximization.

% 各セクションの概要
The rest of the paper is organized as follows. In Section \ref{sec:fram}, we introduce a concept of data-driven MFTD for the simultaneous optimization of material distributions and HF modeling parameters using the multi-channel image, focusing on the significant shift of image data architecture. In Section \ref{sec:for}, we describe the brief formulation of the optimization problem discussed in this study, including the LF optimization and HF evaluation. In Section \ref{sec:num}, we elaborate the detailed implementation of the proposed framework, including the architecture of the multi-channel VAE and the EA-based optimization process. In Section \ref{sec:res}, we examine the effectivenss of the proposed framework through numerical examples of stiffness maximization and maximum stress minimization. Finally, we conclude the paper in Section \ref{sec:con}.

%sssssssssssssssssssssssssssssssssssssssssssssssssssssssssssssssssssssssssssssss
\section{Framework}
\label{sec:fram}
% multifidelity topology design
% 最適化プロセス概要
% DDM、Smoothing、二値化はプリプロセスとして纏める？詳細はAppendixでもOK
%sssssssssssssssssssssssssssssssssssssssssssssssssssssssssssssssssssssssssssssss

To facilitate the simultaneous optimization of the material distributions and HF modeling parameters, we first introduce a concept of MFTD specifically defined for the HF modeling parameter search.
A multi-objective topology optimization problem is formulated in general, as follows:

\begin{equation}
\begin{aligned}
\makebox[2cm][l]{$\underset{\gamma}{\text{minimize}}$}   & [J_1(\gamma), J_2(\gamma), ..., J_{r_\text{o}}(\gamma)] \\
\makebox[2cm][l]{subject to} & G_j(\gamma) \leq 0, \quad (j = 1, 2, ..., r_c), \\
                             & \gamma(\boldsymbol{x}) = 0 \ \text{or} \ 1, \quad \forall \boldsymbol {x} \in D,
\end{aligned}
\end{equation}
\noindent where $J_i(\gamma)$ and $G_j(\gamma)$ represent objective and constraint functionals, respectively.
$\gamma$ denotes the design variable field, which is a binary material density field determined by the spatial posiion $\boldsymbol{x}$ in the design domain $D$, where $D$ denotes a predefined two- or three-dimensional design domain.
The design variable field $\gamma$ comprises binary material density field,
In MFTD, solving this topology optimization problem directly is assumed to be difficult or even impossible in some cases, as it involves dealing with an enormous number of design variables and is formulated as a non-linear mathematical optimization problem.

\begin{figure}[H]
    \centering
    \includegraphics[width=0.7\columnwidth]{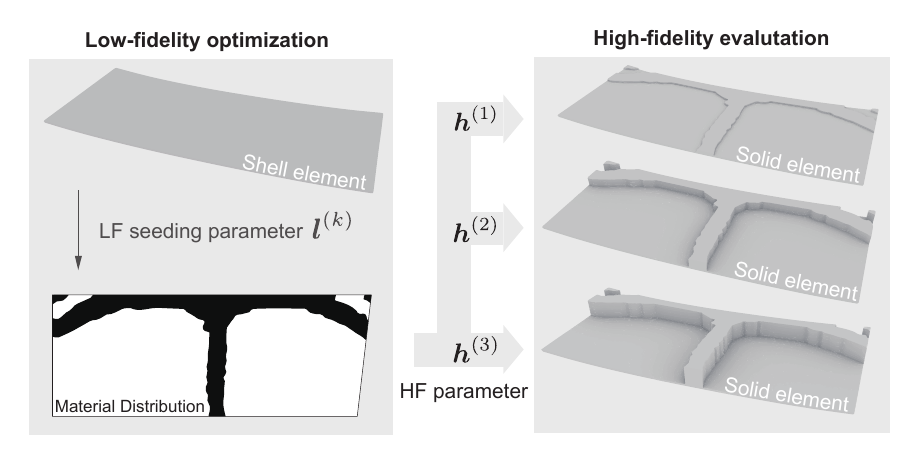}
    \caption{Illustration of LF and HF parameter concept.}
    \label{fig:parameter_concept}
\end{figure}

To indirectly tackle the original intricate topology optimization problem, MFTD divides the original optimization problem into two components: a LF topology optimization and a HF evaluation.
The LF topology optimization problem is formulated, as follows:

\begin{equation}
    \begin{aligned}
    \makebox[2cm][l]{$\underset{\gamma}{\text{minimize}}$}   & [\Tilde{J}_1(\gamma^{(k)}, \boldsymbol {l}^{(k)}), \Tilde{J}_2(\gamma^{(k)}, \boldsymbol {l}^{(k)}), ..., \Tilde{J}_{r_\text{o}}(\gamma^{(k)}, \boldsymbol {l}^{(k)})] \\
    \makebox[2cm][l]{subject to} & \tilde{G}_j(\gamma^{(k)}, \boldsymbol {l}^{(k)}) \leq 0, \quad (j = 1, 2, ..., r_c), \\
                                & 0 \leq \gamma^{(k)}(\boldsymbol{x}) \leq 1, \quad \forall \boldsymbol{x} \in D, \\
    \makebox[2cm][l]{for given}  & \boldsymbol {l}^{(k)},
    \end{aligned}
\label{eq:low-fidelity}
\end{equation}
\noindent
where $\Tilde{J}_i$ and $\Tilde{G}_j$ represent the objective and constraint functionals for the LF optimization problem, respectively.
The vector $\boldsymbol{l} = [l_1, l_2, ..., l_{N_{\text{lf}}^{\text{sd}}}]^\text{T}$ consists of LF seeding parameters, which are used to generate a variety of design candidates by varying the LF topology optimization problem.
The superscript $k$ with $k = 1, 2, ..., N_{\text{lf}}^{\text{sd}}$ corresponds to different topology optimizations with respect to the points in $\boldsymbol{l}$.
In other words, the optimization problem (\ref{eq:low-fidelity}) is solved a total of $N_{\text{lf}}^{\text{sd}}$ times for each time using a different sample point from $\boldsymbol{l}^{(k)}$, resulting in $N_{\text{lf}}^{\text{sd}}$ distinct topologically optimized material distributions $\gamma_{\text{to}}^{(k)}$.
Once all candidates from the LF topology optimization problems are collected, the HF evaluation is performed to identify the best solution, denoted as $\gamma_{\text{to}}^{(k^*)}$, where:
\begin{equation}
k^* \in \underset{k}{\text{arg min}} \{ [J_1(\gamma_{\text{to}}^{(k)}, \boldsymbol{h}^{(k)}), J_2(\gamma_{\text{to}}^{(k)}, \boldsymbol{h}^{(k)} ), ..., J_{r_\text{o}}(\gamma_{\text{to}}^{(k)}, \boldsymbol{h}^{(k)} )] \ | \ G_j (\gamma_{\text{to}}^{(k)}, \boldsymbol{h}^{(k)}) \leq 0\}.
\label{eq:high-fidelity}
\end{equation}
\noindent
Herein, HF seeding parameters are introduced as an extension to the conventional MFTD.
These parameters are represented as a list of vectors denoted by $\boldsymbol{h} = {[ \boldsymbol{h}_1, \boldsymbol{h}_2, ..., \boldsymbol{h}_{N_{\text{hf}}^{\text{sd}}} ]}^\text{T}$.
Each vector $\boldsymbol{h}_m$ contains parameter distributions used to configure the Finite Element (FE) models for HF evaluation.
The HF seeding parameters can include various configurations, such as thickness distributions for creating three-dimensional FE models by extruding the material distribution in the horizontal plane to a specified height in the vertical direction, as illustrated in Fig.~\ref{fig:parameter_concept}.
The original objective and constraint functionals, $J_i$ and $G_j$, are utilized to find the best solutions in the HF evaluation.
When $i \geq 2$, $J_i$ exhibits trade-offs with each other, resulting in the best solutions being Pareto optima represented by the vector $k^*$ that contains multiple sample point indices.

% どうせならLF変数も複数パターンに変えた図の方が良いかも

% low- forward vs high- inverse の利点
%In MFTD, it is crucial that the HF evaluation is exclusively performed through forward analysis, such as static analysis, while LF topology optimization is conducted using iterative inverse analysis.
%This approach simplifies the original optimization problem, which is often challenging to tackle directly with inverse analysis, whether from a theoretical or computational standpoint.
%Forward analysis, on the other hand, offers a more straightforward solution.

% LF & HF seeding parameter両方の探索が重要
The effectiveness of achieving a satisfactory solution set depends on the diversity of candidates collected, considering a range of joint seeding parameters $(\boldsymbol{l},\boldsymbol{h})$.
Notably, both LF and HF seeding parameters have an impact on structural performance and should be optimized to eventually obtain high-performance structures.
The conventional data-driven MFTD focuses on searching for satisfactory solutions in the material distribution space, that is initially determined by LF seeding parameters and automatically updated by the EA-based optimization process, as HF models are evaluated with a consistent configuration applied across all candidates.
On the contrary, our proposed framework updates the initial candidates derived from various sets of LF and HF seeding parameters utilizing the EA-based optimization process, similar to the conventional one but specifically adaptive to the data architecture of a multi-channel image.

\begin{figure}[hbt]
    \centering
    \includegraphics[width=0.7\columnwidth]{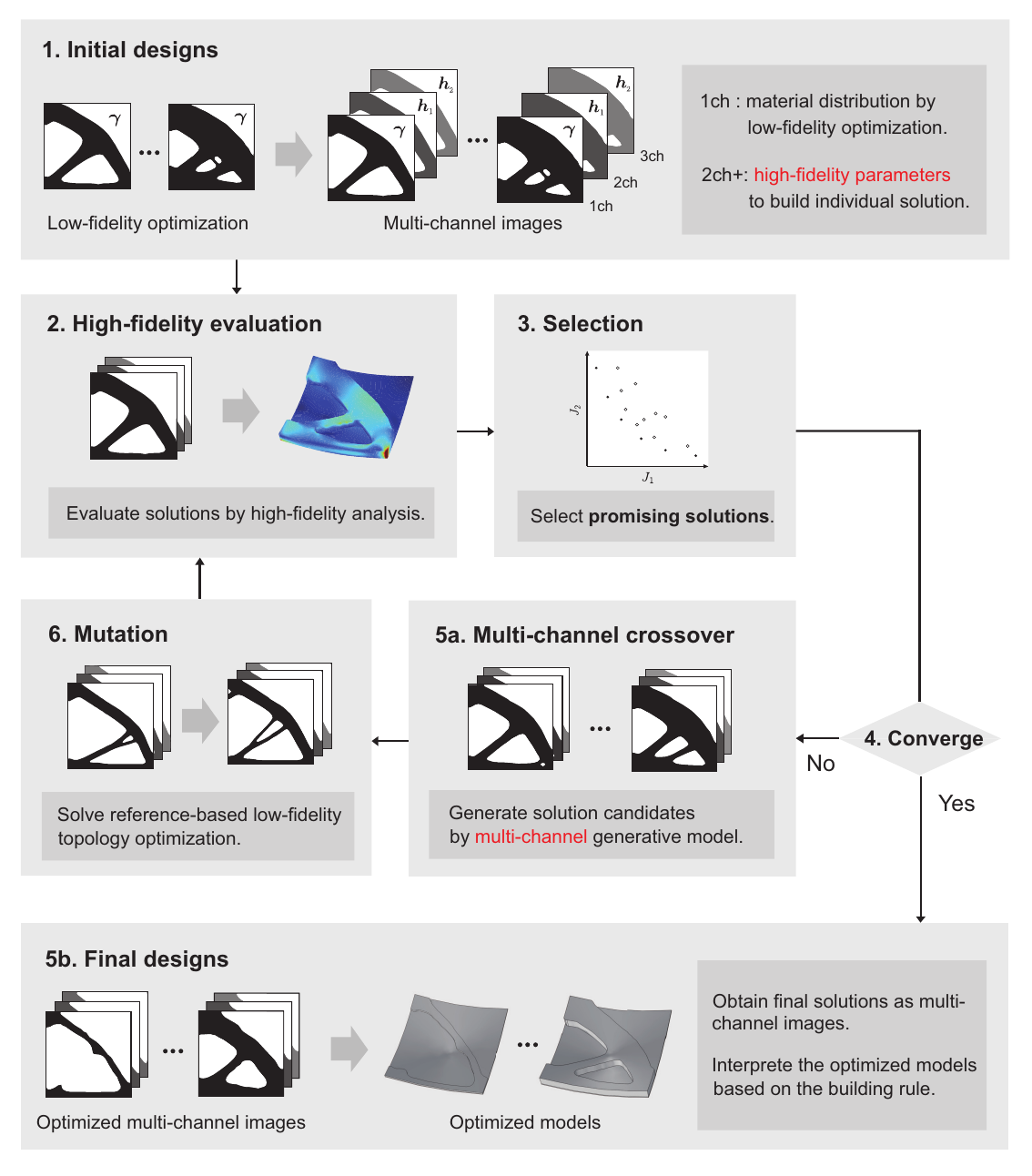}
    \caption{Schematic flowchart of the proposed data-driven multifidelity topology design.}
    \label{fig:flowchart}
\end{figure}

% Overview of proposed framework: initial designs
The overview of our proposed framework is illustrated in Fig.~\ref{fig:flowchart}.
Initially, a dataset of multi-channel images is prepared, encapsulating the necessary information for FE model construction, such as material distributions and thickness of reinforcements.
This preparation involves solving the LF topology optimization (\ref{eq:low-fidelity}) across a range of LF parameters to generate a variety of material distribution candidates, which are then stored in the first channel of the images.
Subsequent channels are dedicated to HF parameters that define the mechanical performance of each solution, which may be represented as either constant values or distributions within the design domain.

% Overview of proposed framework: HF evaluation
Following the initial setup, steps 2 through 6 are executed in a loop until a convergence criterion is met.
Step 2 involves evaluating each candidate solution through forward analysis on the HF model (e.g., static analysis in this context) to derive precise objective functionals.
Here, material distributions from the first channel are mapped onto a three-dimensional design domain, with the HF model being constructed based on thickness values from the second channel.
%Unlike the conventional D2MFTD, which was limited to two-dimensional domains, this study extends the mapping to three-dimensional spaces, accommodating structures that require HF space exploration, such as thickness optimizations, typically modeled with thin, three-dimensional shell elements. (DDMに言及しない場合はコメントアウト)
It should be noted that this framework has the flexibility to optimize not only thickness but also other physical properties through various HF model constructions.

% Overview of proposed framework: EA-based iterative process
Steps 3 to 6 encompass EA-based methods for autonomously refining design candidates.
Step 3 selects promising solutions via elitism, while step 4 assesses convergence using the hypervolume indicator, a metric widely employed in multi-objective optimization to gauge the quality of optimized solution.
If convergence is not achieved, the process advances to step 5a, where new candidates are generated through a generative model in a crossover-like operation, or to step 5b to finalize designs upon convergence.
Whereas the conventional data-driven MFTD utilized a single-channel VAE (SC-VAE) for candidate generation, our framework employs a multi-channel VAE (MC-VAE) \cite{antelmi2019sparse} to enable the neural network of VAE to learn the interactive relationship between material distributions and HF parameters.
Periodically, the solution dataset undergoes mutation-like operation to introduce mutants that significantly diverge from the existing dataset yet likely survive the selection process.
This mutation is guided by an additional constraint referencing the current dataset, with outcomes recorded in the first channel of the mutated multi-images.
Although this discussion focuses on bi-objective scenarios, the framework is adaptable to problems with three or more objectives.

% convergence criterion
As for step 4, similar to the data-driven MFTD based on the SC-VAE, a convergence of the objective space is checked based on a hypervolume indicator for the selected solution set $\gamma_{\text{to}}^{(k^*)}$ using a reference point $r_{\text{hv}} \in \mathbb{R}^{r_o}$, defined as follows \cite{shang2020survey}:
\begin{equation}
\text{HV}(k^*, r_{\text{hv}}) = \mathcal{L} \left(\bigcup_{k_p \in k^*} \{ J(\gamma_{\text{to}}^{(k)}) \ | \ J(\gamma_{\text{to}}^{(k_p)}) \leq J(\gamma_{\text{to}}^{(k)}) \leq r_{\text{hv}} \} \right)
\label{eq:convergence}
\end{equation}
\noindent
where the Lebesgue measure, denoted by $\mathcal{L}(\cdot)$, is utilized to quantify the hypervolume indicator.
Specifically to bi-objective problems, this indicator measures the area spanned by the non-dominated solutions, identified by $k_p$, and a reference point $r_{\text{hv}}$.
The optimization process is deemed to have converged when the hypervolume indicator $\text{HV}$ falls below a predefined threshold $\varepsilon_{\text{HV}}$.

%===========================================================================
\subsection{Selection}
\label{subsec:sel}
% 少し冗長な気がする（特に新しいことは述べていない）
%===========================================================================

As introduced in the previous data-driven MFTD research \cite{yamasaki2021data}, the selection of promising solutions is performed using a concept of elitism through the application of the non-dominated sorting genetic algorithm II (NSGA-II) \cite{deb2002fast}, which is a well-established algorithm in the realm of multi-objective EAs.
NSGA-II focuses on generating a set of non-dominated solutions, satisfying the following relationship, especially in a multi-objective minimization problem:
\begin{equation}
\forall i : J_i (\gamma_{\text{to}}^{(k_1)}) \leq J_i (\gamma_{\text{to}}^{(k_2)}) \quad \land \quad \exists j : J_j (\gamma_{\text{to}}^{(k_1)}) < J_j (\gamma_{\text{to}}^{(k_2)}).
\label{eq:nsga2}
\end{equation}
\noindent
This relationship signifies that a solution $\gamma_{\text{to}}^{(k_1)}$ is not dominated by another solution $\gamma_{\text{to}}^{(k_2)}$, meaning the $k_1$-th solution is considered superior to the $k_2$-th one.

NSGA-II selects solutions using non-dominated sorting and crowding distance sorting.
In the non-dominated sorting phase, a set of solutions is categorized into different fronts based on the dominance relationship (\ref{eq:nsga2}).
The most dominant solutions, i.e., those not dominated by any other solution, are assigned to the first front called {\it Pareto front}.
Successive fronts contain solutions that are dominated only by those in the previous fronts. 
Once solutions are sorted into fronts, the crowding distance for each solution is calculated.
This metric measures the density of solutions surrounding a particular solution in the objective space.
Within each front, solutions with larger crowding distances are given higher priority, as they are more isolated from other solutions and represent more diverse parts of the objective space.
Finally, when selecting solutions for the next generation, NSGA-II starts with those in the lowest rank (first front) and proceeds to higher ranks.
Within a rank, solutions with higher crowding distances are preferred.
This process continues until the desired number of solutions for the next generation is achieved.

%===========================================================================
\subsection{Crossover}
\label{subsec:crossover}
% 
%===========================================================================

%===========================================================================
\subsubsection{Preprocessing image data}
\label{subsec:preprocess}
% 
%===========================================================================

Before employing the MC-VAE to generate a new candidate solution set, the input data containing material distribution information undergo preprocessing to enable a deep generative model to efficiently learn the characteristics of the input data. This preprocessing consists of two steps: Design Domain Mapping (DDM) and the smoothing of the material distributions.

% DDM：3Dメッシュに適用したのは新規性？
DDM, as proposed by Yamasaki et al., serves to map material distributions from the original domain to another domain while preserving the boundary conditions inherited from the original domain \cite{yamasaki2018knowledge}.
In the previous data-driven MFTD, DDM was employed to resize the original data to fit within a unit square domain denoted as $\bar{D}$, which was discretized with square elements, i.e., structured quad mesh, for VAE's convenience.
While the previous research focused primarily on resizing the mesh to a unit square, we extend its focus to mapping a 3D mesh onto a 2D plane mesh, in addition to mesh resizing.
The detailed methodologies behind this dimensional extension is elaborated in the appendix \ref{ap:ddm}.

After mapping the material distributions from the original design domain to the unit plane mesh, these mapped material distributions undergo smoothing via the Helmholtz PDE \cite{lazarov2011filters}, as formulated below:
\begin{equation}
    -r^2 \nabla^2 \bar{\gamma} + \bar{\gamma} = \gamma.
\label{eq:helmholtz}
\end{equation}
\noindent
Here, $\gamma$ represents the mapped material distribution, $\bar{\gamma}$ is the smoothed material distribution, and $r$ denotes the filtering radius.
The smoothed material distribution $\bar{\gamma}$ is computed by solving Eq.~(\ref{eq:helmholtz}) using FEM.

%===========================================================================
\subsubsection{Multi-channel variational auto-encoder}
\label{subsec:mc_vae}
% 
%===========================================================================

% MC-VAEの特徴
Building upon the fundamental principles of the standard VAE, a MC-VAE adds an extra layer of complexity by simultaneously processing multiple channels of input data.
The primary advantage of the MC-VAE lies in its capacity to capture the underlying relationships between different data channels, which would be challenging for separate SC-VAEs to accomplish.

\begin{figure}[hbt!]
\centering
\includegraphics[width=0.8\textwidth]{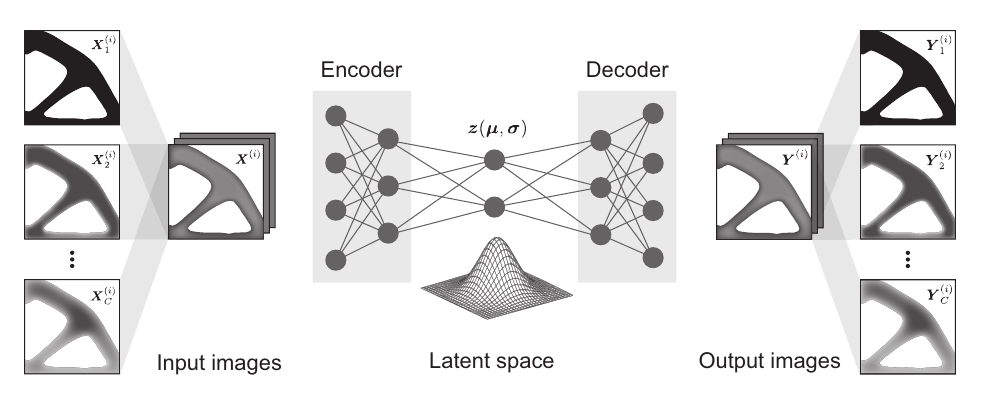}
\caption{Illustration of the data structure of MC-VAE.}
\label{fig:vae_concept}
\end{figure}

% multi-channelの理論
Consider a set of samples $\boldsymbol{X} = \{\boldsymbol{X}^{(1)}, \boldsymbol{X}^{(2)}, \ldots, \boldsymbol{X}^{(N)}\}$ where $\boldsymbol{X}^{(i)} \in \mathbb{R}^{ N_{\text{mc}}}$ denotes normalized pixel data with the total number of pixels $N_{\text{mc}}$, and multiple channels of input data represented as $\boldsymbol{X}^{(i)} = \{ \boldsymbol{X}_1^{(i)}, \boldsymbol{X}_2^{(i)}, ... \boldsymbol{X}_C^{(i)} \}$, where $C$ is the total number of channels, and each $\boldsymbol{X}_c^{(i)} \in \mathbb{R}^{N_{\text{sc}}}$ corresponds to the data from the $c$-th channel with each channel's number of pixels $N_{\text{sc}}$.
The total number of pixels $N_{\text{mc}}$ accordingly results in $C N_{\text{sc}}$.
The first channel $\boldsymbol{X}_1^{(i)}$ contains the material distribution for the sample $\boldsymbol{X}^{(i)}$, while the second and subsequent channels $\boldsymbol{X}_{c+1}^{(i)}$ contain each HF seeding parameters' vector $\boldsymbol{h}_c$.
Fig.~\ref{fig:vae_concept} depicts the concept of the MC-VAE's data structure.

VAEs are generally structured around two principal neural networks: an encoder that transforms high-dimensional input data $\boldsymbol{X}^{(i)}$ into latent variables $\boldsymbol{z} \in \mathbb{R}^{N_{\text{lt}}}$, where $N_{\text{lt}}$ is significantly smaller than $N_{\text{mc}}$, and a decoder that reconstructs high-dimensional output data $\boldsymbol{Y}^{(i)} \in \mathbb{R}^{N_{\text{mc}}}$ from these latent variables.
Within the MC-VAE framework, data from each channel is processed by its dedicated encoder network, yielding several sets of latent variables, namely ${ \boldsymbol{z}_1, \boldsymbol{z}_2, ... \boldsymbol{z}_C }$.
These sets are then concatenated to form a unified latent representation $\boldsymbol{z}$.
The MC-VAE's decoder network is engineered to utilize this unified latent representation to reconstruct the data for each channel, producing $ \boldsymbol{Y}^{(i)} = \{ \boldsymbol{Y}_1^{(i)}, \boldsymbol{Y}_2^{(i)}, ... \boldsymbol{Y}_C^{(i)} \}$.
Following encoding and decoding, the first channel of the decoded image is remapped back onto the original domain to represent the material distribution.
The VAE generates new data through its decoder, facilitated by ensuring that the compressed data adheres to a Gaussian distribution within the latent space.

% 参考文献と照らし合わせて再確認
The fundamental idea of the VAE lies in finding a probability distribution $P(\boldsymbol{X})$ from an approximated model $P_{\boldsymbol{\theta}} (\boldsymbol{X}) = P(\boldsymbol{X}; \boldsymbol{\theta})$, where $\boldsymbol{\theta}$ are the parameters of the decoder to be optimized through learning process so that the likelihood of $P_{\boldsymbol{\theta}} (\boldsymbol{X})$ is maximized ultimately.
In other words, the likelihood to generate data points similar to the input data $\boldsymbol{X}$ is desired to be maximized during the iterations.
However, it is hard to calculate the likelihood $P_{\boldsymbol{\theta}} (\boldsymbol{X})$ due to computational efficiency, so VAEs employ an approach known as variational inference, specifically leveraging the `variational lower bound.'
This approach lets us optimize the expression on the right side of the following inequality instead of maximizing the $P_{\boldsymbol{\theta}} (\boldsymbol{X})$ directly:
\begin{equation}
    \log P_{\boldsymbol{\theta}} (\boldsymbol{X}^{(i)}) \geq - D_{\text{KL}} (Q_{\boldsymbol{\phi}}(\boldsymbol{z}|\boldsymbol{X}^{(i)}) || P_{\boldsymbol{\theta}} (\boldsymbol{z})) + \mathbb{E}_{Q_{\boldsymbol{\phi}}(\boldsymbol{z}|\boldsymbol{X}^{(i)})} \left[\log P_{\boldsymbol{\theta}} (\boldsymbol{X}^{(i)} | \boldsymbol{z})\right],
    \label{eq:variational_lower_bound}
\end{equation}
\noindent
where $Q_{\boldsymbol{\phi}}(\boldsymbol{z}|\boldsymbol{X}^{(i)})$ and $P_{\boldsymbol{\theta}} (\boldsymbol{X}^{(i)} | \boldsymbol{z})$ respectively correspond to the encoder and decoder model where $\boldsymbol{\phi}$ and $\boldsymbol{\theta}$ are the parameters of the encoder and decoder, respectively.
The latent variable, denoted by $\boldsymbol{z}$, is given by:
\begin{equation}
    \boldsymbol{z}^{(i,j)} = \boldsymbol{\mu}^{(i)} + \boldsymbol{\sigma}^{(i)} \odot \boldsymbol{\epsilon}^{(j)},
    \label{eq:reparameterization_trick}
\end{equation}
where the mean and standard deviation of the approximate posterior, denoted as $\boldsymbol{\mu}^{(i)}$ and $\boldsymbol{\sigma}^{(i)}$, respectively, represent the encoder's outputs for each individual sample.
The symbol $\odot$ signifies element-wise multiplication, and $\epsilon^{(j)}$ denotes a random sample drawn from the standard normal distribution.
Eq.~(\ref{eq:reparameterization_trick}) is commonly referred to as the 'reparameterization trick,' a technique that facilitates gradient computation in probabilistic sampling within a neural network.
It accomplishes this by breaking down the stochastic variable into two components: a deterministic element (such as the output of a neural network) and a sample originating from an independent noise source.
It is important to note that the latent space captures features relevant to all channels while still preserving the characteristics of the traditional VAE.

The first term, $D_{\text{KL}}(Q_{\boldsymbol{\phi}}(\boldsymbol{z}|\boldsymbol{X}^{(i)}) || P_{\boldsymbol{\theta}} (\boldsymbol{z}))$, on the right-hand side of Eq.~(\ref{eq:variational_lower_bound}), represents the Kullback–Leibler (KL) divergence, used to quantify the disparity between the encoder's distribution and a standard Gaussian distribution.
Minimizing this term leads to an approximation of the encoder distributions that closely resemble Gaussian distributions. The KL divergence is calculated according to its definition, as follows:

\begin{equation}
\begin{aligned}
    D_{\text{KL}} (Q_{\boldsymbol{\phi}}(\boldsymbol{z}|\boldsymbol{X}^{(i)}) || P_{\boldsymbol{\theta}} (\boldsymbol{z})) & = - \frac{1}{2} \sum_{j=1}^{N_{\text{lt}}} \left(1 + \log((\sigma^{(i)}_j)^2) - (\mu^{(i)}_j)^2 - (\sigma^{(i)}_j)^2\right).
\end{aligned}
\label{eq:KL_divergence}
\end{equation}
\noindent
The second term in Eq.~(\ref{eq:variational_lower_bound}) represents the reconstruction error between the original input material distributions and the decoded probability distributions, computed as binary cross-entropy for each channel, as follows:
\begin{equation}
\begin{aligned}
    \mathbb{E}_{Q_{\boldsymbol{\phi}}(\boldsymbol{z}|\boldsymbol{X}_c^{(i)})} \left[\log P_{\boldsymbol{\theta}} (\boldsymbol{X}^{(i)} | \boldsymbol{z})\right] &= \frac{1}{N_{\text{lt}}} \sum_{j=1}^{N_{\text{lt}}} \log P_{\boldsymbol{\theta}} (\boldsymbol{X}_c^{(i)} | \boldsymbol{z}^{(i,j)}) \\
    & = \frac{1}{N_{\text{lt}}} \sum_{j=1}^{N_{\text{lt}}} \left( \boldsymbol{X}_c^{(i)} \log (\boldsymbol{Y}_c^{(i,j)}) + (\boldsymbol{\text{I}}-\boldsymbol{X}_c^{(i)}) \log(\boldsymbol{\text{I}}-\boldsymbol{Y}_c^{(i,j)}) \right),
\end{aligned}
\label{eq:reconstruction_error}
\end{equation}
\noindent
where $\boldsymbol{\text{I}}$ denotes the identity matrix.
Finally, a loss function $L_{\text{VAE}}$ is computed by introducing a parameter to adjust the weight of the KL divergence, as shown below:
\begin{equation}
    L_{\text{VAE}}(\boldsymbol{\theta}, \boldsymbol{\phi}; \boldsymbol{X}^{(i)}) = w_{\text{KL}} \cdot D_{\text{KL}} (Q_{\boldsymbol{\phi}}(\boldsymbol{z}|\boldsymbol{X}^{(i)}) || P_{\boldsymbol{\theta}} (\boldsymbol{z})) - 
\sum_{c=1}^{C} \mathbb{E}_{Q_{\boldsymbol{\phi}}(\boldsymbol{z}|\boldsymbol{X}_c^{(i)})} \left[\log P_{\boldsymbol{\theta}} (\boldsymbol{X}_c^{(i)} | \boldsymbol{z})\right],
\label{eq:loss_vae}
\end{equation}
\noindent
where $w_{\text{KL}} \in \mathbb{R}^{+}$ represents the weight parameter known as KL weight, which governs the influence of the KL divergence, encouraging the distribution in the latent space to approximate a standard normal distribution.
A higher KL weight urges the model to prioritize aligning the latent space distribution closer to the standard normal distribution.
The value of the KL weight impacts the trade-off between generative diversity and reconstruction capability of the model.
As shown in Eq.~(\ref{eq:loss_vae}), the reconstruction error is computed for each channel first, and then the errors from all the channels are summed.

%===========================================================================
\subsection{Mutation}
\label{subsec:mutation1}
% 
%===========================================================================

% HF seeding parameterありの突然変異について説明する
To prevent premature convergence and maintain dataset diversity throughout the entire design process, mutant material distributions for the first channel are generated by a mutation-like operation which Yaji et al. found effective in promoting more promising candidates.
The desired mutation-like operation involves generating new material distributions that inherit topological features from the material distributions of reference samples in the target solution set under optimization, so that the mutant material distributions can survive the selection process more likely.
This operation is achieved by solving the LF topology optimization problem defined in Eq.~(\ref{eq:low-fidelity}), incorporating an additional constraint function:

\begin{equation}
    \Tilde{G}_{\text{mut}}(\gamma^{(m)}) = \int_D (1-|\gamma^{(m)}-\gamma_{\text{ref}}^{(m)}|) \ d\Omega \ \leq \ \Tilde{G}_{\text{mut}}^{\text{max}} \int_D d\Omega,
\label{eq:mutation_constraint}
\end{equation}
\noindent
where $m=1,2,...,N_{\text{mut}}$, and $\Tilde{G}_{\text{mut}}^{\text{max}}$ is a parameter controlling the discrepancy between a selected reference material distribution $\gamma_{\text{ref}}^{(m)}$ and the optimized material distribution for mutation.
The concept behind this mutation process is detailed in the reference \cite{yaji2022data}.
Notably, a smaller value of $\Tilde{G}_{\text{mut}}^{\text{max}}$ results in a more divergent material distribution, which inherits fewer topological features from the original material distribution.
Careful consideration is required when determining the value of this parameter, as a highly divergent material distribution faces greater challenges in surviving during the selection process.

%sssssssssssssssssssssssssssssssssssssssssssssssssssssssssssssssssssssssssssssss
\section{Formulation}
\label{sec:for}
% トポロジー最適化の定式化
%sssssssssssssssssssssssssssssssssssssssssssssssssssssssssssssssssssssssssssssss

%===========================================================================
\subsection{Topology optimization}
\label{subsec:top}
% トポロジー最適化の基礎理論
% 密度法の考え方：連続変数の導入、フィルタ、プロジェクションの導入
% チェーンルールだけ表示して、随伴法とMMAの詳細には踏み込まない。感度解析には随伴法を使用する旨を紹介
%===========================================================================

This section briefly explore the formulation of topology optimization, specific to the solid isotropic material with penalty (SIMP) method, one of the most widely used methods in topology optimization.
To solve the optimization problem such as (\ref{eq:low-fidelity}), the design variable field $\gamma$ is discretized into a set of design variables $\gamma_e \in [0, 1]$ for each element $e$, where $D \subseteq \mathbb{R}^{d}$ with $d \in {2,3}$ denotes a predefined design domain, using Finite Element Method (FEM) in this study.
The governing equations are formulated based on the SIMP method, detailed as follows:
\begin{equation}
\begin{aligned}
\boldsymbol{K}\boldsymbol{U} & = \boldsymbol{F}, \\
\boldsymbol{K} & = \sum_e^{N_\text{elm} } E_{\text{SIMP}} (\hat{\Tilde{\gamma}}_e) K_e, \\
E_{\text{SIMP}} (\hat{\Tilde{\gamma}}_e) & = E_{\text{min}} + \hat{\Tilde{\gamma}}_e^p (E_0 - E_{\text{min}}),
\label{eq:stiffness_equation}
\end{aligned}
\end{equation}
\noindent
where the stiffness matrix $\boldsymbol{K}$ represents a system of stiffness equations for $N_\text{elm}$ elements, with a displacement vector, denoted by $\boldsymbol{U}$, containing the displacements at each node or element, subject to applied loads represented by a force vector $F$.
Here, $K_e$ denotes an element stiffness matrix, and $E_{\text{SIMP}}$ represents the modified Young's modulus.
$E_0$ and $E_{\text{min}}$ correspond to the Young's modulus of solid and void materials, respectively.
The parameter $p$ serves as a penalty factor that encourages the binarization of design variables.
Furthermore, $\hat{\Tilde{\gamma}}_e$ denotes a filtered and projected design variable, as elaborated below.

To mitigate the emergence of overly complex topologies during optimization, a filtering process is typically incorporated into topology optimization.
The design variables are subjected to filtering using a method proposed by Bourdin \cite{bourdin2001filters} to smooth the design variables by averaging each one with its neighbors within a filtering radius $r_{\text{f}}$, as shown below:
\begin{equation}
\begin{aligned}
   \Tilde{\gamma}_e & = \frac{\sum_{J \in \mathbb{N}_{d,e}} H(\boldsymbol {x}_j)\gamma_j}{\sum_{J \in \mathbb{N}_{d,e}} H(\boldsymbol {x}_j)}, \\
H(\boldsymbol {x}_j) & = \text{max}(0, r_{\text{f}} - |\boldsymbol {x}_j - \boldsymbol {x}_e| ),
\end{aligned}
\label{eq:filter_design_variable}
\end{equation}
\noindent
where $\mathbb{N}_{d,e}$ represents a set of design variables located within a distance less than the filter radius $r_{\text{f}}$ from the position of each element $\boldsymbol {x}_e$.
Additionally, a Heaviside projection is incorporated to remove the grayscale results and achieve simpler yet sharper topologies \cite{wang2011projection}, as illustrated below:
\begin{equation}
\begin{aligned}
   \hat{\Tilde{\gamma}}_e & = \frac{\text{tanh}(\beta (\Tilde{\gamma}_e - \eta)) + \text{tanh}(\beta \eta)}{\text{tanh}(\beta (1 - \eta)) + \text{tanh}(\beta \eta)},
\end{aligned}
\label{eq:project_design_variable}
\end{equation}
\noindent
where $\beta$ represents a scaling parameter that determines the sharpness of the projection, and $\eta$ is a threshold value to determine whether the projected variable should be adjusted closer to zero or one.
The design variable field $\gamma$ is iteratively updated using the Method of Moving Asymptotes (MMA), a well-established technique of sequential convex programming method \cite{svanberg1987method}, based on the sensitivity of the objective functional against the design variable field, calculated using the chain rule. 

%===========================================================================
\subsection{Stiffness maximization}
\label{subsec:dar}
% LF最適化
% HF評価
%===========================================================================

A bi-objective topology optimization to maximize stiffness and minimize volume is written as follows:

\begin{equation}
    \begin{aligned}
    \makebox[2cm][r]{find} \quad   & \gamma_e \quad (e = 1,2,...,N_\text{elm} ), \\
    \makebox[2cm][r]{that minimize} \quad & W = \boldsymbol{U}^\text{T} \boldsymbol{K} \boldsymbol{U},\\
    \quad & V = \sum_{e=1}^{N_\text{elm} } v_e \gamma_e, \\
    \makebox[2cm][r]{subject to} \quad & \gamma_e \in \{ 0, 1\},
    \end{aligned}
\label{eq:bi_stiffness}
\end{equation}
\noindent
where $v_e$ represents the volume of each element.
The objective functional $W$ is defined as the compliance of the structure, or strain energy, calculated by $\boldsymbol{U}^\text{T} \boldsymbol{K} \boldsymbol{U}$, which can be considered as the inverse of stiffness.
Notably, the material distribution derived from this formulation is entirely discrete by its nature, allowing for a more precise evaluation of objective functionals than methods involving relaxation.

% LF optimization
The bi-objective optimization problem (\ref{eq:bi_stiffness}) is converted to a mono-objective optimization problem with a constraint of maximum volume $V_{\text{max}}$, as follows:

\begin{equation}
    \begin{aligned}
    \makebox[2cm][r]{find} \quad   & \gamma_e \quad (e = 1,2,...,N_\text{elm} ), \\
    \makebox[2cm][r]{that minimize} \quad & W = \boldsymbol{U}^\text{T} \boldsymbol{K} \boldsymbol{U},\\
    \makebox[2cm][r]{subject to} \quad & V = \sum_{e=1}^{N_\text{elm} } v_e \gamma_e \ \leq \ V_{\text{max}}, \\
    \quad & \gamma_e \in [ 0, 1 ],
    \end{aligned}
\label{eq:mono_stiffness}
\end{equation}
\noindent
where $\gamma_e$ is transitioned from a discrete to a continuous variable.
Solving this optimization problem (\ref{eq:mono_stiffness}) with various values of the maximum volume constraint $V_{\text{max}}$ allows us to generate a variety of material distribution candidates to be stored in the first channel of multi-channel images as the initial designs, combined with randomly assigned HF seeding parameters.

% HF evaluation
The material distributions derived from LF topology optimization or the first channel of the MC-VAE output undergo a process of smoothing to approximate more binary-like distributions, which are subsequently binarized using a density threshold $\gamma_e = 0.5$. 
This binarization process ensures that the material distribution is clearly defined as either solid or void, with only the elements designated as solid being extruded to the height of the HF seeding parameters, stored in the second channel of the multi-channel images, to design the reinforcement structures.

%===========================================================================
\subsection{Maximum stress minimization}
\label{subsec:opt}
% LF最適化
% HF評価
%===========================================================================

A bi-objective topology optimization to minimize maximum stress and volume is written as follows:

\begin{equation}
    \begin{aligned}
    \makebox[2cm][r]{find} \quad   & \gamma_e \quad (e = 1,2,...,N_\text{elm}), \\
    \makebox[2cm][r]{that minimize} \quad & \sigma_{\text{max}} = \text{max}({\sigma}_{\text{vm},e}),\\
    \quad & V = \sum_{e=1}^{N_\text{elm} } v_e \gamma_e, \\
    \makebox[2cm][r]{subject to} \quad & \gamma_e \in \{ 0, 1\},
    \end{aligned}
\label{eq:bi_stress}
\end{equation}
\noindent
where $\sigma_{\text{vm},e}$ is von Mises stress of each element. 

% LF optimization
This bi-objective optimization problem (\ref{eq:bi_stress}) is so-called a mini-max problem, which is not applicable to gradient-based optimization methods due to its non-differentiability.
To allow for gradient-based optimization, a $p$-norm function \cite{yang1996stress, duysinx1998new} is employed to approximate the maximum stress within the design domain, transforming the original bi-objective optimization problem (\ref{eq:bi_stress}) into a mono-objective optimization problem, as follows:

\begin{equation}
    \begin{aligned}
    \makebox[2cm][r]{find} \quad   & \gamma_e \quad (e = 1,2,...,N_\text{elm} ), \\
    \makebox[2cm][r]{that minimize} \quad & \sigma_{\text{PN}} = \left( \frac{1}{N_\text{elm} } \sum_{e=1}^{N_\text{elm} } \hat{\sigma}_{\text{vm},e}^P (\gamma_e) \right)^{\frac{1}{P}},\\
    \makebox[2cm][r]{subject to} \quad & V = \sum_{e=1}^{N_\text{elm} } v_e \gamma_e \ \leq \ V_{\text{max}} \ , \\
    %\quad & \hat{\sigma}_{\text{vm},e} (\gamma_e) = \gamma_e^q \sigma_{\text{vm},e} \ , \\
    \quad & \gamma_e \in [ 0, 1 ],
    \end{aligned}
\label{eq:mono_stress}
\end{equation}
\noindent
where the stress field $\hat{\sigma}_{\text{vm},e} (\gamma_e) = \gamma_e^q \sigma_{\text{vm},e}$ is modified to facilitate the binarization of design variables and mitigate numerical challenges encountered when optimizing design variables near material boundaries, where stress can change abruptly.
This modification ensures that design variables transition smoothly, avoiding sharp increases as their values approach $0$.
The approximated stress, denoted by $\sigma_{\text{PN}}$, converges to the actual maximum stress $\sigma_{\text{max}}$ as $P$ approaches infinity.
However, as $P$ increases, so does the computational complexity.
Consequently, the value of $P$ should be carefully determined compromising the computational stability and evaluation accuracy.

% HF evaluation
Similar to the approach for stiffness maximization, material distributions are binarized, and HF models are constructed using the defined material boundaries and HF parameters.
Subsequently, von Mises stress $\sigma_{\text{max}}$ and volume $V$ are evaluated within these FE models.
It is important to note that the maximum stress is determined directly from the von Mises stress distribution across the design domain, bypassing any relaxation or approximation strategies.

%sssssssssssssssssssssssssssssssssssssssssssssssssssssssssssssssssssssssssssssss
\section{Numerical implementation}
\label{sec:num}
% 補強部材の設計コンセプトについて説明：外板上に描画された断面形状を板厚方向に押し出して設計
% アルゴリズムの擬似コード
% Too muchでなければ設計プロセスごとの説明。サブセクション設けない分はセクション下に纏めて記述
%sssssssssssssssssssssssssssssssssssssssssssssssssssssssssssssssssssssssssssssss

This section outlines the detailed implementation of the proposed framework.
Building upon the core concepts previously discussed, a pseudo-code of the proposed framework is presented in Algorithm (\ref{alg:overall}), tailored for bi-objective optimization problems such as those defined by Eq.~(\ref{eq:bi_stiffness}) and Eq.~(\ref{eq:bi_stress}).
The parameters governing the entire process are detailed in Table \ref{tab:overall}.
The LF seeding parameter $\boldsymbol{l}$ and the HF seeding parameter $\boldsymbol{h}$ are determined through Latin Hypercube Sampling (LHS) from the intervals $[l^{\text{min}}, l^{\text{max}}]$ and $[h^{\text{min}}, h^{\text{max}}]$, respectively.
LHS is a statistical technique designed to generate a diverse set of plausible parameter combinations from a multidimensional distribution, commonly employed in computer simulations to efficiently explore varied design spaces. %LHS
The bi-objective functionals $J_1$ and $J_2$ are evaluated under the condition that the constraint functional $G$ is met for the dataset.
The design iteration continues until the relative error of the hypervolume indicator falls below $\varepsilon_{\text{HV}}$, or the iteration count reaches a predefined maximum $N_{\text{max}}$.
The reference point for hypervolume calculations, as per Eq.~(\ref{eq:convergence}), is set to values marginally exceeding the highest initial objective values found on the Pareto front.
A mutation-like process is conducted at intervals of $N_{\text{mut}}^{\text{int}}$ iterations, solving the LF topology optimization for $N_{\text{mut}}$ reference samples across $N_{\text{mut}}^{\text{sd}}$ variations of LF seeding parameters.
Regarding the mutation constraint (\ref{eq:mutation_constraint}), the reference material distributions $\gamma_{\text{ref}}^{(m)} (m=1,2,...N_{\text{mut}})$ are uniformly selected from the Pareto front solutions denoted as $\Theta^{(i)}$.
Among all the processes, HF evaluation demands the most computational resources, as it involves running numerical simulations for each sample in the temporary dataset $\Theta_{\text{tmp}}$ at every iteration $i$.
To optimize computation time, HF evaluations (lines 9-11) were performed in parallel, leveraging the independence of each simulation from the others.

% selection
In the selection process, the number of offspring should be carefully determined due to its influence on the efficiency of both the entire optimization and MC-VAE's learning performance.
The MC-VAE's learning performance is significantly influenced by both the network architecture and the volume of the input dataset, particularly when this dataset is relatively limited.
A consistent offspring count for the MC-VAE's input dataset promotes more stable learning outcomes because the MC-VAE's network architecture remains unchanged throughout the optimization iterations.
On the contrary, the more offspring count allow the MC-VAE to learn the characteristics of the input dataset more effectively.
Therefore, we opt to determine the number of offspring by setting a minimum threshold and selecting all Pareto front points as offspring when their total number surpasses this minimum.

% Software列挙
The framework outlined in Algorithm (\ref{alg:overall}) was developed using Python (version 3.9.13) as the primary programming environment, incorporating a variety of tools for its implementation.
The LF topology optimization and HF evaluation, along with the solution of Helmholtz PDEs, were executed in COMSOL Multiphysics (version 6.1), a commercial FEM software.
In this study, COMSOL Multiphysics was controlled via scripts written in MATLAB (version 2023a).
The MC-VAE was constructed using TensorFlow (version 2.12.0), an open-source machine learning framework renowned for its versatility in developing and deploying sophisticated machine learning models.
DDM was accomplished with the aid of a mapping function from the IGL Python library (version 2.4.1).
The traditional EA operations, such as selection and hypervolume calculations, were facilitated by DEAP (version 1.3.3), an open-source Python library designed for crafting and executing genetic algorithms and other evolutionary computing strategies.

\begin{algorithm}[hbt!]
    \caption{Pseudo-code of the expanded data-driven MFTD for a bi-objective optimization problem}
    \begin{algorithmic}[1]
        \For{$k = 1$ to $N_{\text{lf}}^{\text{sd}}$}
            \State Solve the LF optimization problem (\ref{eq:mono_stiffness}) or (\ref{eq:mono_stress}) to find $\gamma^{(k)}_e (e=1, 2, ..., )$ on $\boldsymbol{l}^{(k)}$
        \EndFor
        \State Contain the material distributions in the 1ch of a temporary dataset, $\Theta_{\text{tmp},1} \leftarrow \{\gamma^{(1)}_e, \gamma^{(2)}_e,..., \gamma^{(N_{\text{lf}}^{\text{sd}})}_e \}$
        \For{$c=2$ to $C$}
            \State Define HF parameters for the $c$-ch of the temporary dataset, $\Theta_{\text{tmp},c} \leftarrow \{\boldsymbol{h}^{(1)}, \boldsymbol{h}^{(2)},..., \boldsymbol{h}^{(N_{\text{lf}}^{\text{sd}})} \}$
        \EndFor
        \For{$i=0$ to $N^{\text{max}}$}
            \For{$k=1$ to $\text{len}(\Theta_{\text{tmp}})$}
                \State Calculate HF evaluation functionals, $\{ J_1^{(i)} (\gamma^{(k)}_e), J_2^{(i)} (\gamma^{(k)}_e),G^{(i)} (\gamma^{(k)}_e) \}$
            \EndFor
        
            \State Squeeze $\Theta_{\text{tmp}}$ so that all the material distributions satisfy $G^{(i)} \leq 0$
            \If{$i=0$}
                \State Set the current dataset, $\Theta^{(i)} \leftarrow \Theta_{\text{tmp}}$
            \Else
                \State Augment the dataset, $\Theta^{(i)} \leftarrow \Theta^{(i)} \cup \Theta_{\text{tmp}}$
            \EndIf
            \State Squeeze $\Theta^{(i)}$ by the selection algorithm based on $\mathcal{J}(\Theta^{(i)}) \leftarrow \bigcup_{k=1}^{\text{len}(\Theta^{(i)})} \{ J_1^{(i)} (\gamma^{(k)}_e), J_2^{(i)} (\gamma^{(k)}_e) \}$
            \If{the hypervolume indicator $\text{HV} \leq \varepsilon_{\text{HV}}$ regarding $\mathcal{J}(\Theta^{(i)})$}
                \State Return the obtained dataset $\Theta^{(i)}$ with $\mathcal{J}(\Theta^{(i)})$
            \EndIf
            \State Convert $\Theta^{(i)}$ to normalized dataset by design domain mapping, $\boldsymbol{X} \leftarrow \{ \boldsymbol{X}^{(1)}, \boldsymbol{X}^{(2)}, ..., \boldsymbol{X}^{(\text{len}(\Theta^{(i)}))} \}$
            \State Acquire a new dataset by MC-VAE using $\boldsymbol{X}$ as the input, $\boldsymbol{X}_{\text{VAE}} \leftarrow \{ \boldsymbol{X}_{\text{VAE}}^{(1)}, \boldsymbol{X}_{\text{VAE}}^{(2)}, ..., \boldsymbol{X}_{\text{VAE}}^{(N_{\text{VAE}})} \}$
            \State Convert $\boldsymbol{X}_{\text{VAE}}$ to $\Theta_{\text{tmp}}$ by the inverse design domain mapping
            \If{$\text{mod}(i,N_{\text{mut}}^{\text{int}})$}
                \State Sample reference material distributions for mutation, $\{ \gamma_{e, \text{ref}}^{(1)}, \gamma_{e, \text{ref}}^{(2)},..., \gamma_{e, \text{ref}}^{(N_{\text{mut}})} \} \in \Theta^{(i)}$
                \For{$n=1$ to $N_{\text{mut}}^{\text{sd}}$}
                    \For{$m=1$ to $N_{\text{mut}}$}
                        \State Solve the LF optimization problem (\ref{eq:mono_stiffness}) or (\ref{eq:mono_stress}) for $\gamma_e^{(m)}$ on $\boldsymbol{l}^{(n)}$ with the constraint (\ref{eq:mutation_constraint})
                    \EndFor
                    \State Assemble a dataset of the optimized material distributions, $\theta_{\text{mut},1}^{(n)} \leftarrow \{\gamma_e^{(1)}, \gamma_e^{(2)},..., \gamma_e^{(N_{\text{mut}})} \}$
                    \For{$c=2$ to $C$}
                        \State Define HF parameters for $c$-ch of the mutation dataset, $\theta_{\text{mut},c}^{(n)} \leftarrow \{ \boldsymbol{h}^{(1)}, \boldsymbol{h}^{(2)},...,\boldsymbol{h}^{(N_{\text{mut}})} \}$
                    \EndFor
                \EndFor
                \State Augment the temporary dataset with all the mutation dataset, ${\Theta_{\text{tmp}}} \leftarrow \Theta_{\text{tmp}} \cup \bigcup_{n=1}^{N_{\text{mut}}^{\text{sd}}} \theta_{\text{mut}}^{(n)}$
            \EndIf
        \EndFor
        \State Return the obtained dataset $\Theta^{(i)}$ with $\mathcal{J}(\Theta^{(i)})$
    \end{algorithmic}
    \label{alg:overall}
\end{algorithm}

\begin{table}[hbt]
    \centering
    \footnotesize
    \caption{Parameters for the overall procedures.}
    \begin{tabular}{l*{9}{c}r}
        \hline
        Parameter & Symbol & Value\\
        \hline
        Number of LF seeding parameters & $N_{\text{lf}}^{\text{sd}}$ & $100$ \\
        Number of channels & $C$ & $2$ \\
        Maximum iterations & $N^{\text{max}}$ & $100$ \\
        Number of seeding parameters in mutation & $N_{\text{mut}}^{\text{sd}}$ & $5$ \\
        Number of mutant generation for each seeding parameter & $N_{\text{mut}}$ & $10$ \\
        Interval of mutation & $N_{\text{mut}}^{\text{int}}$ & $5$ \\
        Threshold value of hypervolume indicator & $\varepsilon_{\text{HF}}$ & $1.0 \times 10^{-5}$ \\
        Number of samples generated by MC-VAE & $N_{\text{VAE}}$ & $256$ \\
        \hline
    \end{tabular}
    \label{tab:overall}
\end{table}

%===========================================================================
\subsection{Low-fidelity topology optimization}
\label{subsec:lf2}
%
%===========================================================================

% Discretization
For both the stiffness maximization problem (\ref{eq:mono_stiffness}) and the maximum stress minimization problem (\ref{eq:mono_stress}), the raw design variable field $\gamma_e$, filtered design variable field $\Tilde{\gamma}_e$ and the displacement field $\boldsymbol{U}$ are represented using shell elements that feature five intergal points: one at each corner of the quadrilateral and one at the center.
In contrast, the projected design variable field $\hat{\Tilde{\gamma}}_e$ are represented using $\mathbb{P}^0$ Lagrange quadrilateral finite elements, which are suited for element-wise material distribution.
The parameters for the LF topology optimization are detailed in Table \ref{tab:low_fidelity}.

% MMA, adjoint method
The design variables are iteratively updated using the MMA, with a move limit set to $0.05$ for all LF optimization simulations.
The sensitivity analysis of the objectives relative to the design variables employs a discrete adjoint method installed in COMSOL Multiphysics.
The maximum number of iterations was set to $50$ for stiffness maximization and increased to $100$ for maximum stress minimization to accommodate the latter's greater non-linearity and complexity.
The $P$ in Eq.~(\ref{eq:mono_stress}) was initially set to $8$, and increased as $16, 32, 64$ every $30$ iterations based on the continuation method \cite{li2015volume}.
The penalty parameter for stress relaxation was set as $q = 0.5$.

\begin{table}[H]
    \centering
    \footnotesize
    \caption{Parameters for LF optimization problem.}
    \begin{tabular}{l*{8}{c}r}
        \hline
        Parameter & Symbol & Value\\
        \hline
        Filter radius & $r_{\text{f}}$ & $0.03$ \\
        Scaling factor of projection & $\beta$ & $4$ \\
        Threshold value of projection & $\eta$ & $0.5$ \\
        Young's modules of solid material & $E_0$ & $1.0$ \\
        Young's modules of void material & $E_{\text{min}}$ & $1.0 \times 10^{-9}$ \\
        Number of elements & $N_{\text{elm}}$ & $64 \times 64$ \\
        Penalty factor & $p$ & $3$ \\
        \hline
    \end{tabular}
    \label{tab:low_fidelity}
\end{table}

%===========================================================================
\subsection{Multi-channel variational auto-encoder}
\label{subsec:mcvae2}
% 
%===========================================================================

% CNNsを導入した動機
Here, we explore the detailed implementation of MC-VAE, including its network architecture.
Unlike the approach taken in previous research \cite{yaji2022data}, which utilized a basic VAE architecture comprising a multilayer perceptron with two dense layers, we adopt CNNs, which is engineered specifically for processing and understanding visual data.
Whereas dense layers flatten the input image's pixels and assign independent weights to each pixel, convolutional layers apply a series of localized dot products across small sections of the input image using a filter.
Furthermore, while the number of weights in a dense layer simply corresponds to the total pixel count of multi-channel images, CNNs significantly reduce the number of weights by utilizing a filter that matches the input image's channel count and sliding this filter across all image pixels.
Consequently, CNNs excel in capturing the spatial relationships within multi-channel data through their convolutional layers.
This ability to extract spatial features enables VAEs to generate more meaningful and compact representations within the latent space, leading to enhanced image reconstruction during decoding.
The powerful feature extraction capabilities of CNNs motivated us to implement MC-VAE using CNNs to more efficiently manage the growing complexity of input data, which now includes multiple channels.

% 複数チャンネルに拡張するにあたって、単チャンネルとの差異や工夫を記す
The encoder of MC-VAE is composed of an input layer accepting $64 \times 64$ pixels and two channels, two convolutional layers with $32$ and $64$ filters, respectively, a flattening layer, and a dense layer.
The decoder consists of an input layer accepting a latent vector of size $16$, a dense layer, and two deconvolutional layers with $64$ and $32$ filters, respectively.
The convolutional and deconvolutional layers have a kernel size of $3 \times 3$ and a stride of $2 \times 2$.
The final layer applies a sigmoid activation function to ensure the output values are between $0$ and $1$, whereas every convolutional layer but the final one applies a ReLu activation function.
It is important that the output values for the first channel can be directly used as material distributions $\gamma \in [0,1]$, but for the second channel, the values needs denormalization that inversely normalizes the output values from $[0,1]$ to $[h^{\text{min}},h^{\text{max}}]$.
As the optimizer in training neural networks, we used Adam \cite{kingma2014adam}, an extension of the stochastic gradient descent method, which is widely used as the optimizer of neural networks including VAEs.
The maximum epochs, the learning rate, and the mini-batch size are set as $300$, $1.0\times 10^{-4}$, and $16$, respectively.
The KL weight $w_{\text{KL}}$ for controlling the influence of the KL divergence is set as $w_{\text{KL}} = 1.0 \times 10^{-3}$.
The neural networks are learned until the epoch is reached to the maximum number of epochs or the loss function $L_{\text{VAE}}$ is not improved for 40 epochs in total.

\begin{figure}[hbt!]
    \centering
    \includegraphics[width=0.9\textwidth]{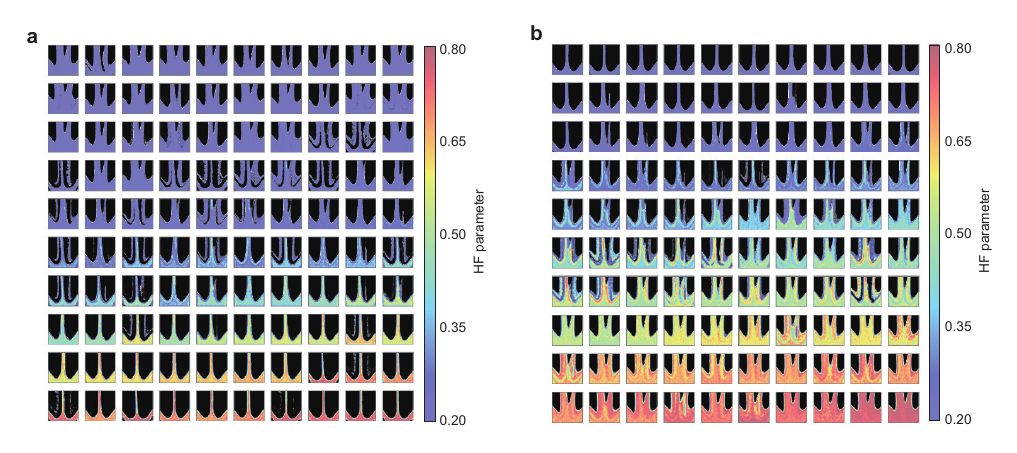}
    \caption{Generated images comparison between a) the proportional and b) inversely proportional dataset.}
    \label{fig:vae_test}
\end{figure}

% volumeと比例・反比例したHFパラメータを与えたときの学習結果を比較する。サンプル数について触れる。一般的な画像学習よりも圧倒的に少なくても学習できているのはサンプルの特徴が非常にシンプルだから
To validate the implementation of the MC-VAE's network architecture, the MC-VAE was trained on two distinct example datasets.
In both datasets, the first channel contained identical material distributions, while the second channel held a scalar HF parameter for each sample.
For one dataset, this scalar HF parameter was directly proportional to the volume fraction derived from the first channel; for the other, it was inversely proportional.
The images generated from these datasets are displayed in Fig.~\ref{fig:vae_test}.
Here, solid materials in the first channel are color-coded to represent the HF parameter in the second channel, with black indicating void areas.
Despite the relatively small size of each dataset, consisting of only $100$ samples —-- significantly less than the typical dataset sizes for VAEs, which are on the order of $10^4$ --— the results, as shown in Fig.~\ref{fig:vae_test}, demonstrate that the MC-VAE successfully captured both the topological features of the material distributions and the relationship between these features and the HF parameter.
This is evident as samples with a higher volume in the proportional dataset correspondingly stored higher HF parameters in the second channel, and the inverse was true for the inversely proportional dataset.
This effective learning from a small dataset can be attributed to the simplicity of the features in the input images and the straightforward relationship between the material distribution and the HF parameter.
Essentially, the MC-VAE needed to grasp either a proportional or inversely proportional relationship between the two channels.
As the number of channels increases, there may be a need to enlarge the input dataset or adjust the network architecture for more efficient learning.
Furthermore, Fig.~\ref{fig:vae_test} suggests that the HF parameters in the second channel are represented not as scalar values but as distributions, even when the input dataset's second channel contained scalar values.
Therefore, when calculating HF models, if a scalar value for the HF parameter is required, it must be extracted from the decoded images as a scalar value through a mathematical computation like averaging.
As a preliminary study, this research utilized scalar values for the HF parameter in each sample to assess the framework's capability to identify promising samples within the HF modeling space.

% イテレーション間で推移するサンプル数について触れる。多少増えて過学習になっても問題ない。
The size of decoding samples is fixed at $256$, while the size of input samples often exceeds $100$, due to the implementation of a minimum offspring number and the inclusion of all Pareto front samples in the selection process.
This excess raises the potential for overfitting within the MC-VAE, as the input sample size surpasses the initial dataset size used for network tuning.
Nonetheless, minor overfitting does not critically impact the framework's effectiveness, as MC-VAE's primary role is to facilitate crossover by producing a multitude of samples that inherit target dataset's features.
These decoded samples, while bearing resemblance to the input dataset in terms of topological features and HF parameters, exhibit slight variations.
It is important to note that the generation of samples with significantly different features is achieved through the mutation-like operation, as detailed in section \ref{subsec:mutation1}.

A particular challenge for MC-VAE is maintaining sample diversity, especially when the input dataset is significantly imbalanced.
The selection of solutions based on LF and HF parameters, coupled with the optimization problem's nature, often results in a disproportionate representation of samples with extreme values for one objective and the inverse for the other.
This imbalance is particularly pronounced for samples with medium values of both objectives, becoming more acute with each iteration and leading to an increasingly biased input dataset for MC-VAE.

% Oversampling
To mitigate the iterative bias in the MC-VAE dataset, we employed a sampling strategy called \textit{oversampling} \cite{ling1998data}, applied at the input stage of MC-VAE.
Oversampling is utilized to rectify class imbalances within datasets by duplicating under-represented samples until a balance is achieved.
This encourages MC-VAE within the data-driven MFTD framework to develop a more comprehensive understanding of less frequent data, thereby enhancing its capability to learn the input dataset's characteristics more effectively maintaining the diversity.

%===========================================================================
\subsection{High-fidelity evaluation}
\label{subsec:hf2}
\begin{figure}[hbt!]
    \centering
    \includegraphics[width=0.6\textwidth]{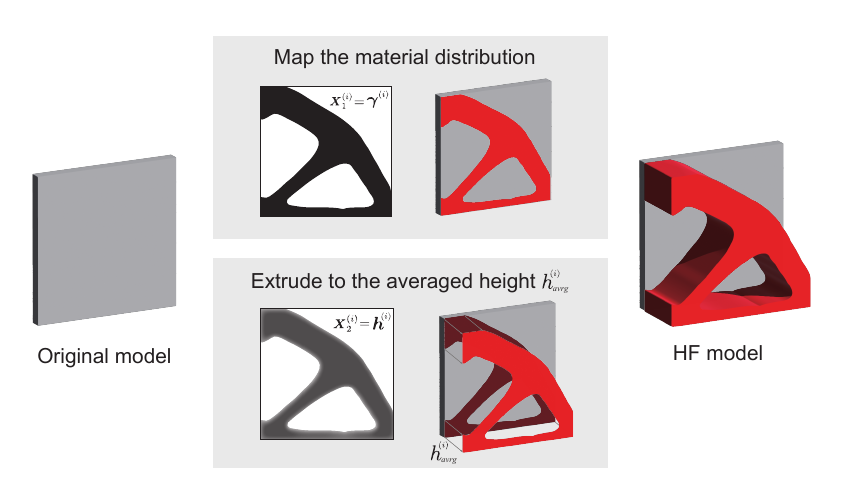}
    \caption{Illustration of the HF modeling.}
    \label{fig:HF_modeling}
\end{figure}

%===========================================================================
% プロセス全体の妥当性を検証するためまずは均一な高さにしたが、本来は分布を扱えるはず。HFモデリングの例を図示したい
In the scenarios of both stiffness maximization and maximum stress minimization, the reinforcement structures are created by extruding the material distributions, thresholded by $\gamma_e \geq 0.5$, to a height determined by the scalar HF parameter value for each sample.
This approach allows each sample's HF parameter field $\boldsymbol{h}^{(k)}$ to contain a scalar value, facilitating the validation of the entire framework, as depicted in Fig.~\ref{fig:HF_modeling}.
It is important to note that the HF evaluation can potentially accommodate HF models with vector HF parameters that vary across the design domain, such as thickness or temperature distributions, owing to the adaptive learning capabilities of MC-VAE.

The displacement and stress fields are discretized using $\mathbb{P}^1$ Lagrange tetrahedral finite elements.
Significantly, the mesh configuration can differ between LF and HF models without issue, leveraging the multifidelity approach where LF optimization and HF evaluation are independently processed.
Consequently, this study employs solid meshes for HF models to ensure detailed evaluations and shell meshes for LF models to conserve computational resources.

%===========================================================================
\subsection{Mutation}
\label{subsec:mutation2}
% hoge
%===========================================================================

\begin{figure}[hbt!]
\centering
\includegraphics[width=0.7\textwidth]{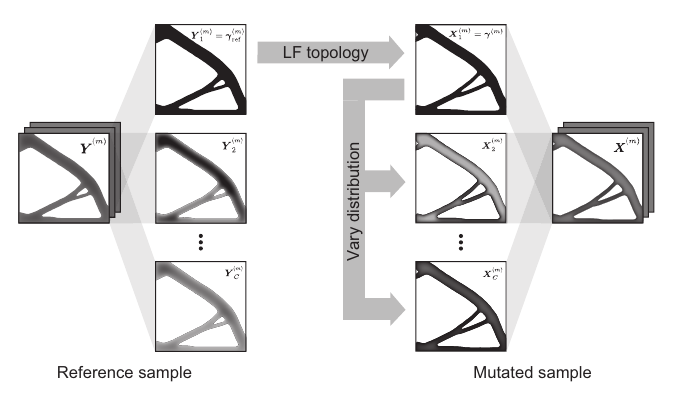}
\caption{Illustration of the mutation process.}
\label{fig:mutation_concept}
\end{figure}

This section explore the mutation-like operation tailored for multi-channel images, which encapsulate both material distributions and HF parameter fields.
Following the generation of diverse multi-channel images via the MC-VAE, a subset of $N_{\text{mut}}$ multi-channel images $\boldsymbol{Y}^{(m)}$ is earmarked for mutation.
The concept of this multi-channel mutation-like operation adapted in this study is depicted in Fig.~\ref{fig:mutation_concept}.
Specifically, the material distribution $\gamma_{e, \text{ref}}^{(m)}$ within the first channel of a reference image $Y_1^{(m)}$ is utilized as a constraint in the LF topology optimization process, as described by Eq.~(\ref{eq:mutation_constraint}), with a discrepancy parameter set to $\Tilde{G}_{\text{mut}}^{\text{max}} = 0.5$.
The rest of the parameters for this optimization mirror those of the standard LF topology optimization.

Upon completion of the LF topology optimization for mutation, the resultant material distribution $\gamma_e^{(m)}$ is contained within the first channel of the mutated sample image as $\boldsymbol{X}_1^{(m)}$.
For subsequent channels in the mutated samples, it is crucial that these distributions correlate with the primary material distribution to facilitate MC-VAE's learning of inter-channel relationships.
Consequently, the HF parameters in these channels $\boldsymbol{X}_c^{(m)}$ are assigned as distributions obtained by multiplying the pixel values in the first channel and the predefined HF parameters.

As outlined in section \ref{subsec:mcvae2}, this study adopted a scalar value for the second channel to assess the framework's efficacy.
Thus, the HF parameter field in the second channel is assigned a distribution derived by multiplying a scalar HF parameter sampled via LHS and the pixel values in the first channel.
To refine this approach, HF parameters in subsequent channels of mutated samples could be directly inherited from the corresponding channels of reference samples, enhancing the probability of these mutated samples surviving through the selection phase.

%sssssssssssssssssssssssssssssssssssssssssssssssssssssssssssssssssssssssssssssss
\section{Results and discussion}
\label{sec:res}
% 数値例により提案手法の有効性を示す
%sssssssssssssssssssssssssssssssssssssssssssssssssssssssssssssssssssssssssssssss

% Stiffness maximization
The efficacy of the proposed framework is illustrated through two numerical examples.
The initial example addresses a stiffness maximization problem, a cornerstone problem within structural design.
This example seeks to demonstrate that the proposed framework can identify high-performance solutions by efficiently exploring the HF parameter space, alongside the traditional LF parameter space that delineates material distributions.
It is presumed that the topology optimization for this straightforward problem can be directly resolved.
For validation, the outcomes derived from the proposed framework are juxtaposed with those obtained from reference solutions solved with uniformly sampled LF and HF parameters.

% Maximum stress minimization
The subsequent example tackles a maximum stress minimization problem, which inherently requires relaxation and approximation methods for stress to be solvable.
Prior research has shown that the mini-max stress problem could be indirectly addressed through the data-driven MFTD, yet without specifically aiming to explore the HF parameter space for superior solutions \cite{kato2023tackling}.
Thus, this second example endeavors to demonstrate the effectiveness of the proposed data-driven MFTD in navigating both LF and HF parameter spaces for the maximum stress minimization problem to secure high-performance solutions.
To confirm this, the stress optimization problem was solved using relaxation and approximation methods, followed by an evaluation of the objective values with HF models, bypassing relaxation and approximation.
Furthermore, as an alternative LF optimization within the data-driven MFTD, the stiffness maximization problem was also solved, subsequently assessing the stress and volume objectives to gauge the impact of the LF optimization problem on the final solution's performance.

%===========================================================================
\subsection{Problem setup}
\label{subsubsec:probset1}
% 問題設定
%===========================================================================

% 設計対象のサーフェス
In practical design scenarios, the target surface is often three-dimensional thin-walled structures with curvature, necessitating the capability of the framework to handle such structures.
To facilitate this, a part of a quadric surface was selected for the target surface, as expressed:
\begin{equation}
f(x,y) = \alpha \sqrt{x^2 + y^2},
\label{eq:cone}
\end{equation}
\noindent
where $\alpha$ was chosen to be $0.3$.

\begin{figure}[hbt!]
    \centering
    \includegraphics[width=1.0\textwidth]{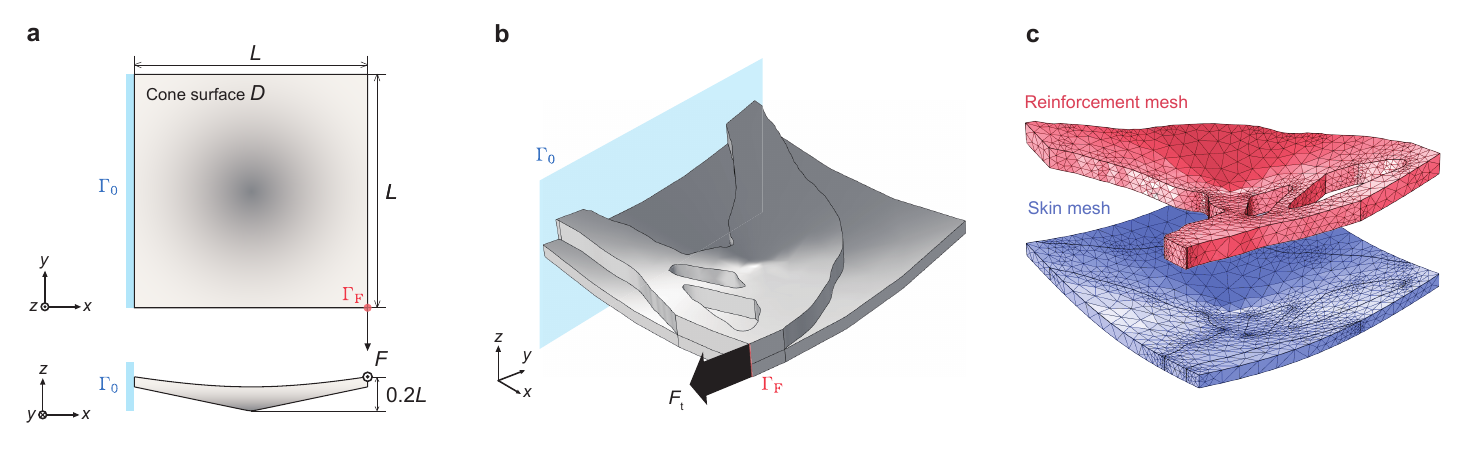}
    \caption{Setup for the LF optimization and HF evaluation for stiffness maximization: a) dimensions and boundary conditions for the LF optimization; b) boundary conditions for a HF evaluation; c) unstructured mesh for the HF evaluation.}
    \label{fig:setup}
\end{figure}

% 各問題設定の図
The setup for LF optimization, including dimensions, loading conditions, and boundary conditions, is depicted in Fig.~\ref{fig:setup} a).
The design domain $D$ is modeled as a cone surface with square-trimmed edges, where the width is defined by the reference length $L$.
A uniform edge load of $F = 1$ is applied at the lower right corner $\Gamma_{\text{F}}$ in the negative $y$ direction.
The Dirichlet boundary condition $\boldsymbol{U} = 0$ is applied to the left wall $\Gamma_{0}$ in the design domain $D$.
This setup is widely recognized as the cantilever beam problem, a prevalent benchmark for topology optimization, albeit with the design domain $D$ being a three-dimensional bending surface.
%It is important to note that while an L-shaped design domain is also a common benchmark for maximum stress minimization, this study exclusively concentrates on the cantilever beam configuration for initial investigations.

The HF configuration, particularly regarding loading and boundary conditions, is depicted in Fig.\ref{fig:setup} b).
These configurations aim to extend the LF optimization settings for thin-walled shell models, as shown in Fig \ref{fig:setup} a), to those suitable for extruded solid models.
Consequently, the edge load in the LF configuration is broadened to an area load in the HF setup.
Tetrahedral meshes for the skin and reinforcement structures are shown in Fig.~\ref{fig:setup} c), colored in blue and red, respectively.
The mesh density is adjusted based on proximity to the boundary surfaces, enhancing the accuracy of stress evaluation in areas where stress is likely to concentrates.

%===========================================================================
\subsection{Validation}
\label{subsec:val}
% 剛性最大化問題
%===========================================================================

%===========================================================================
\subsubsection{Comparison with reference solutions}
\label{subsubsec:comparison1}
% LF結果・初期解と隣り合わせの図に
% Hypervolumeの推移
% 参照解との性能比較：博論より後に解析した参照解を使用
% MC-VAEからデコードされた画像の例
% 最適解の可視化：参照解と最適解から一つずつ
%===========================================================================

\begin{figure}[hbt!]
\centering
\includegraphics[width=.9\textwidth]{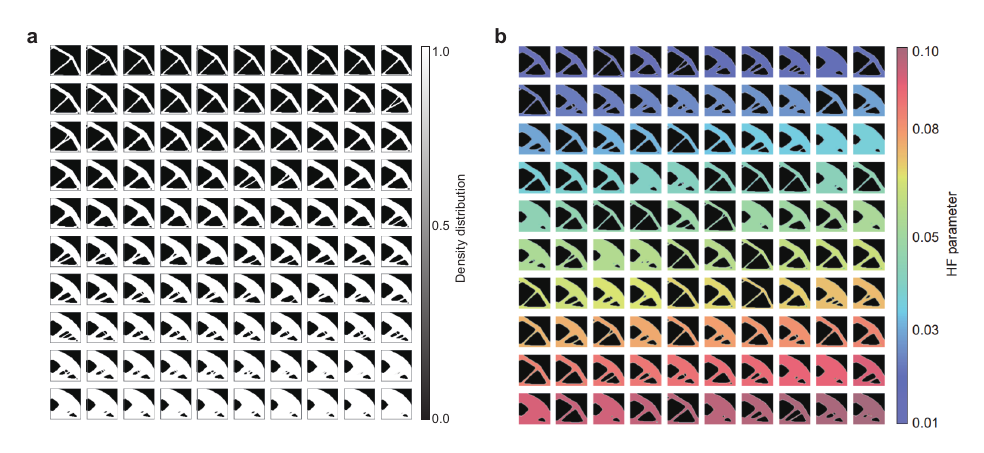}
\caption{Material distributions used as initial designs: a) material distributions obtained from LF stiffness maximization; b) initial designs constructed by assinging HF parameters to corresponding material distributions.}
\label{fig:LF_results_stiffness}
\end{figure}

% LF結果の考察
The material distributions derived from the LF optimization are depicted in Fig.~\ref{fig:LF_results_stiffness} a), organized by the maximum volume constraint $V_{\text{max}}$ ranging from $0.2$ to $0.8$.
Samples with higher $V_{\text{max}}$ exhibit a thicker topology of material distribution.
The central ribs, connecting the arch-shaped outer rib to the bottom rib, are notably thinner compared to the connecting ribs.
These configurations are characteristic outcomes of the cantilever beam problem addressed through topology optimization, specifically using the SIMP method, albeit with slightly thicker and curved ribs at the cone's apex.
This curvature is attributed to the top of the cone being more prone to deformation without significant support, leading to a potential surge in the strain energy, represented by $\boldsymbol{U}^\text{T} \boldsymbol{K} \boldsymbol{U}$.
Consequently, the results from the LF stiffness maximization align qualitatively with the anticipated outcomes.

% 初期解の可視化
The initial dataset was created by initializing the high-fidelity HF parameters using LHS for both $h^{(k)}$ and $l^{(k)} = V_{\text{max}}$, as illustrated in Fig.~\ref{fig:LF_results_stiffness} b).
The samples are organized based on the assigned HF parameters.
Fig.~\ref{fig:LF_results_stiffness} b) demonstrates that the initial dataset comprehensively represents all categories of HF parameters, thanks to the even distribution achieved by LHS.
This dataset was then prepared for input into the MC-VAE, excluding samples that did not meet the $V_{\text{max}}$ constraint, as material boundaries might shift during the smoothing phase.

% LF結果＆初期解

\begin{figure}[hbt!]
\centering
\includegraphics[width=0.5\textwidth]{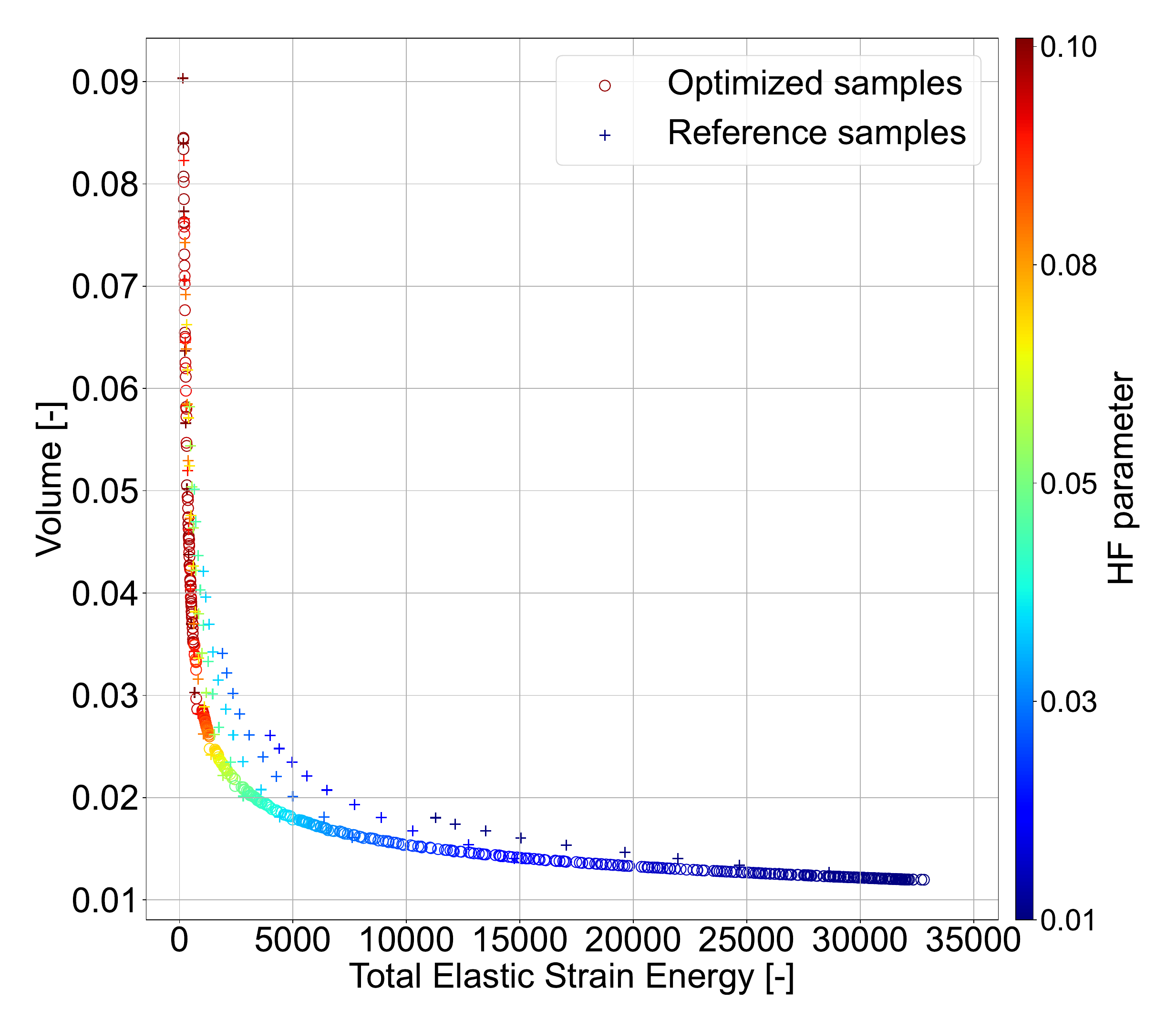}
\caption{Comparison of objective functionals with color of HF parameters between the optimized and reference samples.}
\label{fig:reference_vs_proposed}
\end{figure}

% 参照解との比較
To validate the proposed framework, $100$ reference solutions were obtained by solving the same stiffness maximization problem, and then extruded to heights uniformly sampled from the same range as the proposed framework's HF parameters.
The reference solutions were then evaluated using the HF model, and the objective values were compared with those from the proposed framework.
Objective values derived from the final solutions of the proposed framework and those from the reference solutions are juxtaposed in Fig.~\ref{fig:reference_vs_proposed}.
This figure indicates the effectiveness of the proposed framework in identifying high-performance solutions across broader objective spaces, as the proposed framework's solutions with lower volume exhibit considerably superior performance compared to the reference solutions, with higher volume's performance being equivalent.
This implies that such a parametric study by uniformly sampled HF parameters can yield suboptimal performance without a comprehensive search for more global solutions through the exploration of enormous HF parameters.

\begin{figure}[hbt!]
\centering
\includegraphics[width=0.5\textwidth]{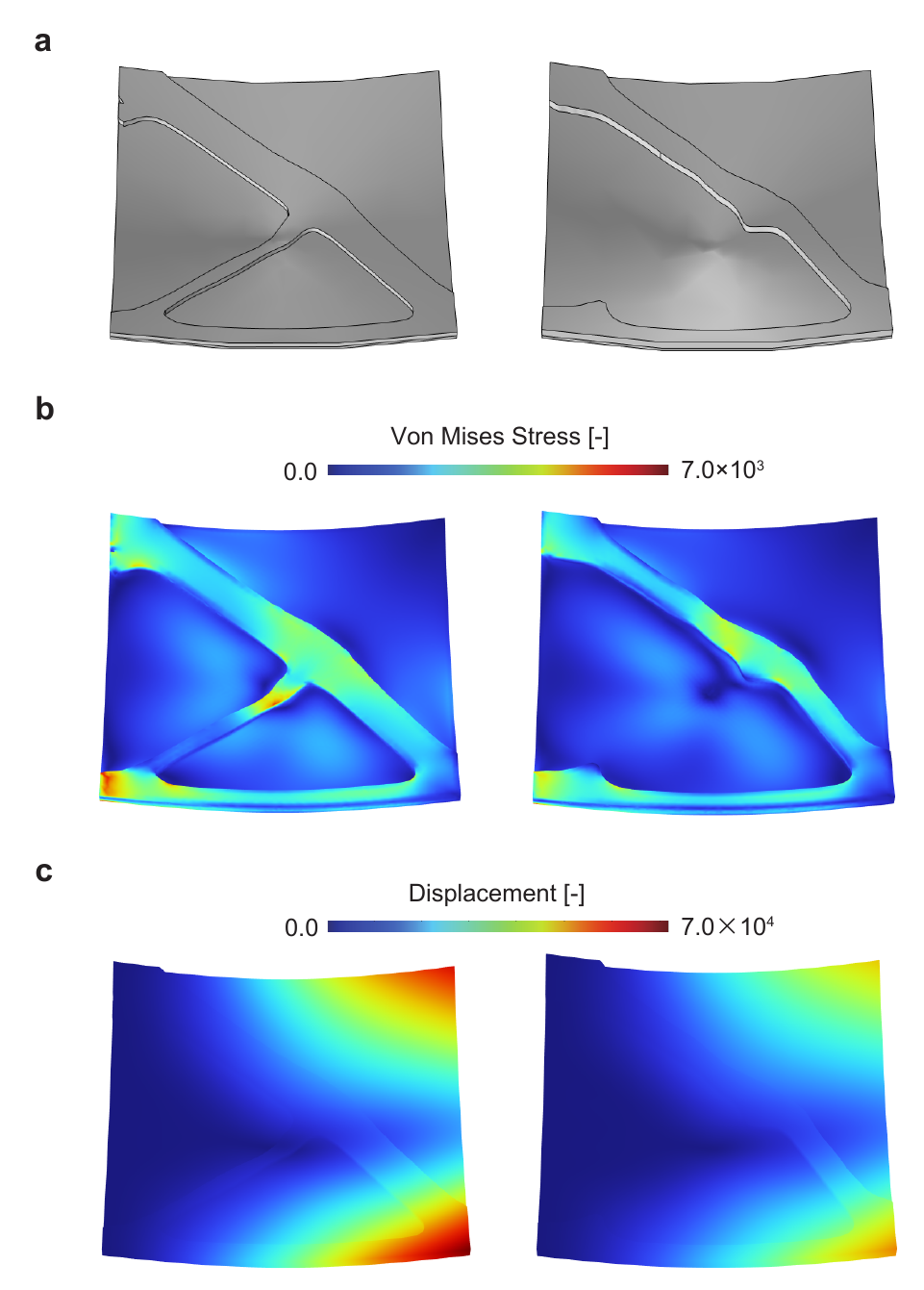}
\caption{Comparison of structural performance between reference (left) and optimized (right) model: a) HF model; b) stress field; c) displacement field. The volume of both models is 0.015, whereas the strain energy is $1.3\times 10^4$ (reference); $9.9\times 10^3$ (optimized).}
\label{fig:reference_vs_proposed_ex}
\end{figure}

% 最適解の可視化
The structural performance of the optimized model is compared with the reference model in Fig.~\ref{fig:reference_vs_proposed_ex} on the right and left, respectively.
Although the volume of both models is same, as shown in Fig.~\ref{fig:reference_vs_proposed_ex} a), the optimized model's reinforcement has a taller and thinner sectional shape without a middle bridge that the reference model has, leading to a $24\%$ lower strain energy than that of the reference model.
This is evident in the stress and displacement fields, as shown in Fig.~\ref{fig:reference_vs_proposed_ex} b) and c), respectively.
Both models exhibit similar displacement and stress distributions, where the displacement is significant at the top and bottom right corner and the stress concentrates on the top and bottom corners of design and the center of the design domain with the most curvature.
This is because their topological features are nearly identical, however, the optimized model's reinforcement is more efficiently supporting the structure with thinner sectional shape and larger height, resulting in the entirely lower stress and displacement fields than the reference model.

%===========================================================================
\subsubsection{Effect of oversampling}
\label{subsubsec:eff_oversampling}
% ハイパーボリューム指標の推移
% MC-VAEからデコードされた画像の例
% 最適解の可視化：参照解と最適解から目的関数の値が近いものを選んで、材料分布・応力分布・変位分布を並べた図
%===========================================================================

\begin{figure}[hbt!]
\centering
\includegraphics[width=0.5\textwidth]{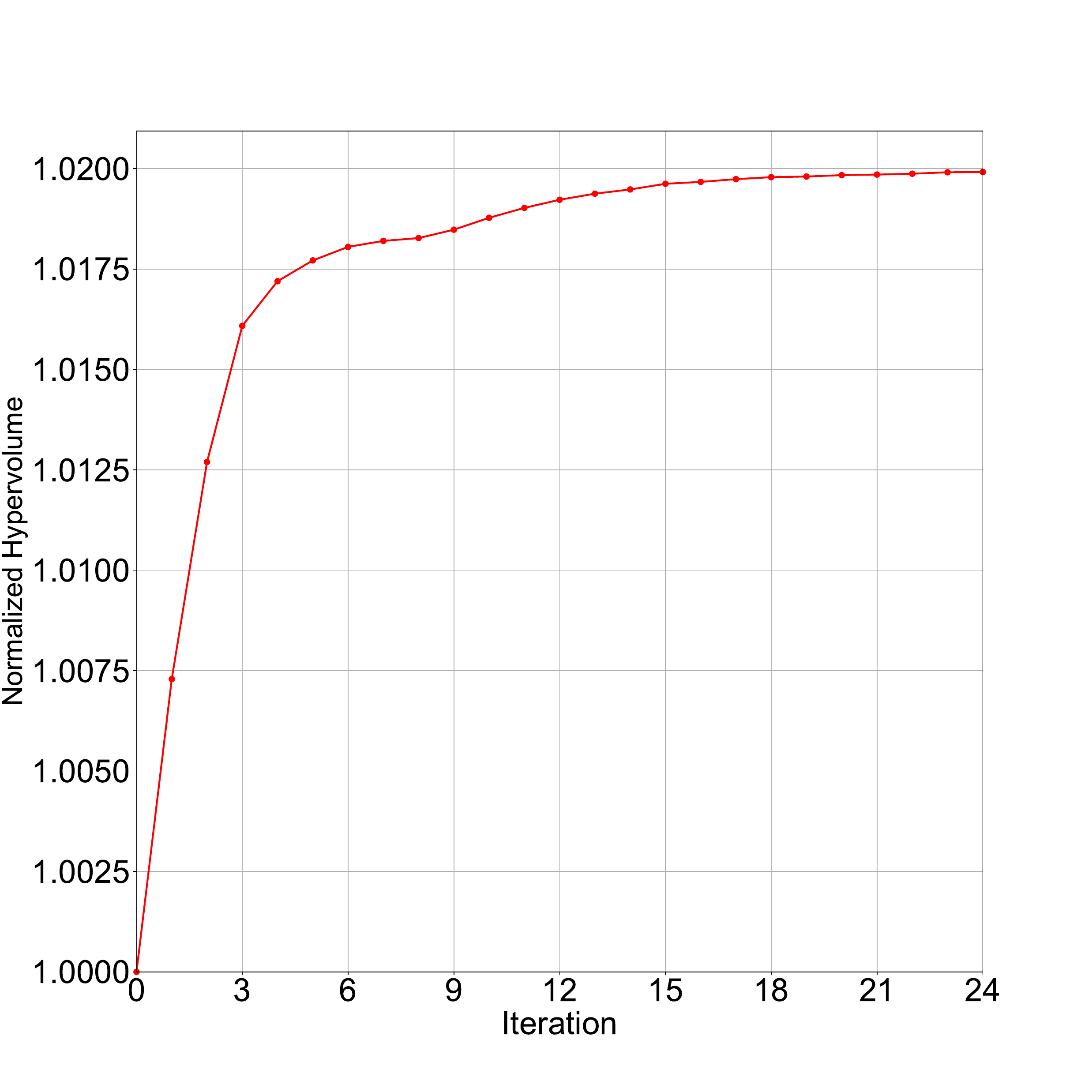}
\caption{Iteration history of hypervolume indicator for stiffness maximization.}
\label{fig:hypervolume_history_stiffness}
\end{figure}

% ハイパーボリュームの遷移
The evolution of the hypervolume indicator from the stiffness maximization is depicted in Fig.~\ref{fig:hypervolume_history_stiffness}.
Hypervolume indicators were computed using $r_{\text{hv}} = [J_1, J_2] = [35000, 0.1]$, where the objective functions $J_1$ and $J_2$ represent strain energy and volume, respectively.
The optimization process converged at the 24th iteration.
This occurred when the relative error of the hypervolume indicator dropped below the threshold.
The hypervolume indicator improved over iterations, indicating an expansion of the Pareto front towards the utopia point in the objective space.
The improvements in the hypervolume indicator were $2.0\%$.
Notably, the convergence of these cases occurred significantly sooner, and the enhancements in the hypervolume indicator were more modest compared to the outcomes for the data-driven MFTD in previous studies \cite{yaji2022data}, which tackled thermal-fluid problems characterized by greater non-linearity than the stiffness maximization problem.
The earlier convergence and lesser improvement in the hypervolume indicator for stiffness maximization should be attributed to the relative simplicity of its optimization problem.

\begin{figure}[hbt!]
\centering
\includegraphics[width=.9\textwidth]{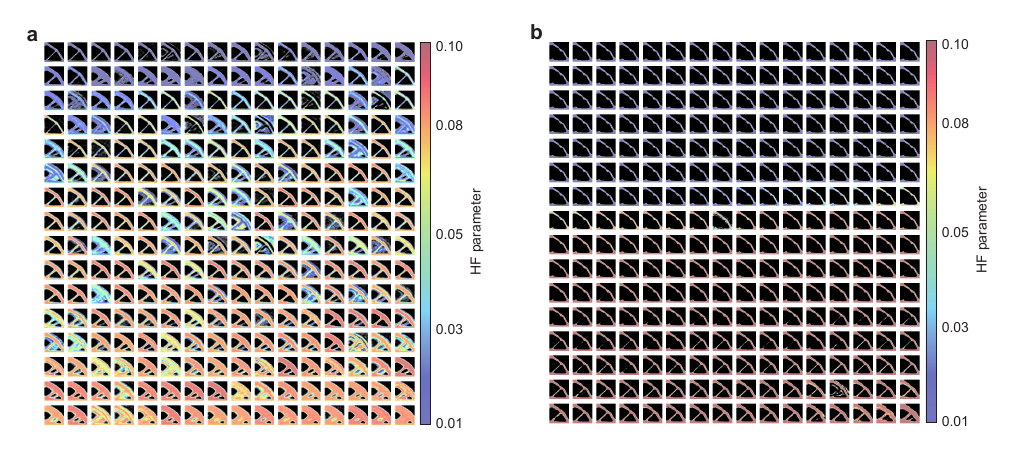}
\caption{Denormalized images decoded by MC-VAE: a) at the initial iteration and b) at the last iteration.}
\label{fig:vae_output}
\end{figure}

% Init and Last VAE 
Fig.\ref{fig:vae_output} displays the denormalized images decoded from the MC-VAE at the a) initial and b) last iteration, where the images are sorted by the corresponding total mass values.
At the last iteration, the simpler topologies with center ribs and no holes are more prevalent, however, as the total mass increases, the topologies feature more complex reinforcements with holes and ribs to support the entire structure more effectively.
This trend is consistent with the height of reinforcement as well; as the total mass increases, the height of reinforcement increases, and the strain energy decreases.
These trends validate that the MC-VAE effectively learns the physical relationship between the sectional shape and height of reinforcements, stored in the first and second channels, respectively.

Additionally, the material boundaries of the material distributions decoded at the last iteration are more distinct than those decoded at the initial iteration.
This discrepancy arises because the dataset for the MC-VAE has no physical relationship between the material distribution and the HF parameter at the initial iteration, increasing the difficulty in effectively learning the relationship between the two channels.

%===========================================================================
\subsection{Maximum stress minimization}
\label{subsec:optimized}
% 最大応力最小化問題
% LF結果：博論から白黒を反転？
% 初期解の可視化
% Hypervolumeの推移
% 従来手法SC-VAEとの性能比較
% LF剛性最大化とLF最大応力最小化との比較
% 最適解の可視化：従来解と提案解から一つずつ
% 潜在変数分布
% LFとHFの目的関数トレンド比較
%===========================================================================

Now we apply the proposed framework, validated through its application to a fundamental topology optimization problem, to the maximum stress minimization that involves strong non-linearity.

%===========================================================================
\subsubsection{Comparison with conventional framework}
\label{subsubsec:comparison2}
%===========================================================================

\begin{figure}[hbt!]
\centering
\includegraphics[width=0.5\textwidth]{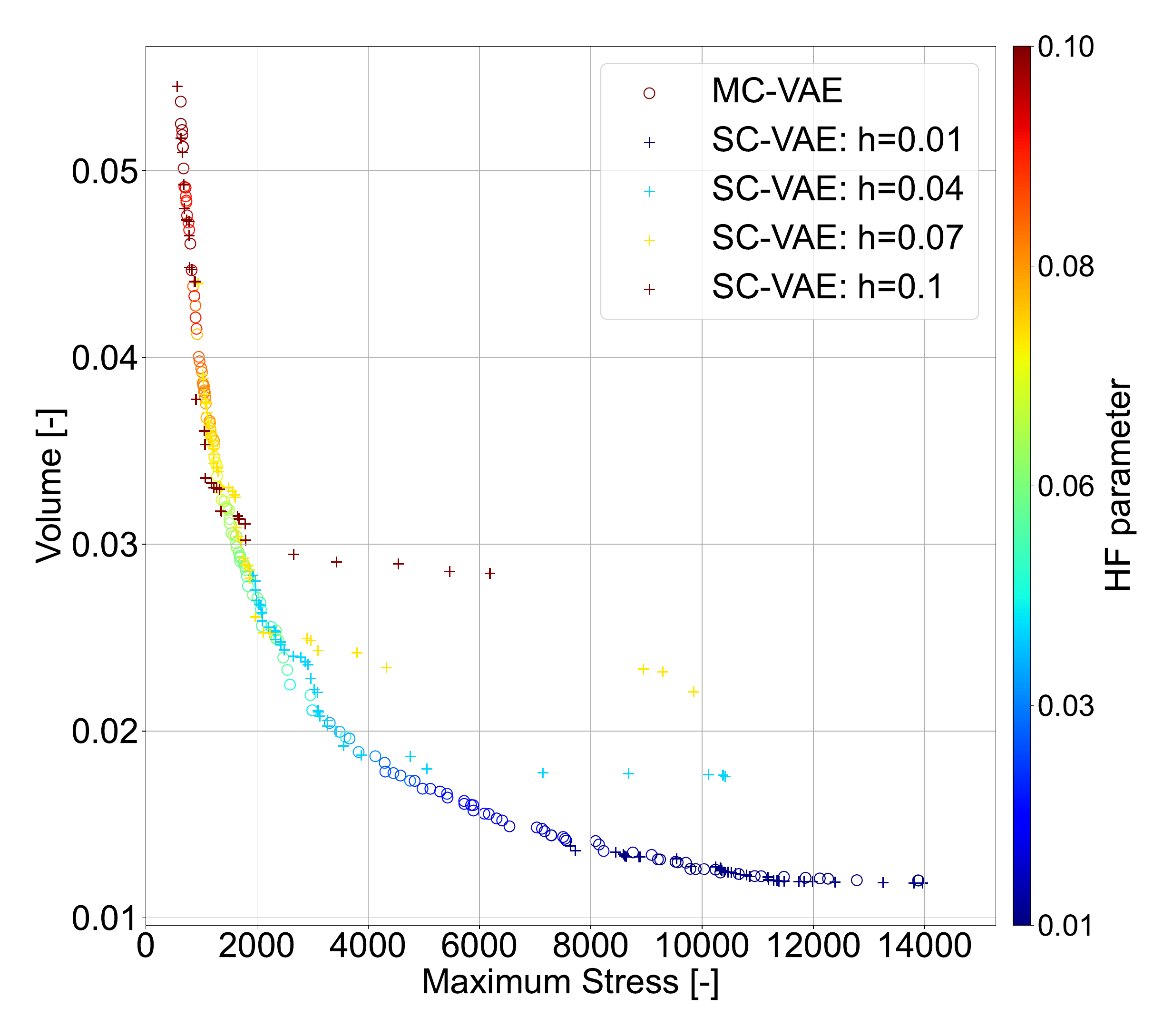}
\caption{Comparison of objective functionals with color of HF parameters between the conventional data-driven MFTD with SC-VAE and the expanded data-driven MFTD with MC-VAE.}
\label{fig:objectives_conventional_vs_proposed}
\end{figure}

To assess the capability of the proposed data-driven MFTD in identifying high-performance samples throughout the HF parameter space, the conventional data-driven MFTD approach was employed for maximum stress minimization using the same problem settings and HF parameters $h=0.01, 0.04, 0.07, 0.1$.
The objective functionals derived from both the conventional data-driven MFTD with SC-VAE and the proposed data-driven MFTD with MC-VAE are illustrated in Fig.~\ref{fig:objectives_conventional_vs_proposed}.
The conventional data-driven MFTD with SC-VAE yielded samples that are partly positioned on the Pareto front for the MC-VAE case, whereas samples with smaller volumes and higher maximum stress are often positioned away from the Pareto front, indicating they are less effective.
Predicting the placement of optimized samples within the objective space prior to applying the conventional data-driven MFTD with a singular HF parameter is challenging.
Therefore, to achieve a Pareto front as comprehensive as that obtained with the proposed data-driven MFTD, extensive parametric studies involving various HF parameters would be necessary if relying solely on the conventional data-driven MFTD.
This underscores the significant advantage of the proposed data-driven MFTD in globally searching for optima across material distributions and HF parameters.
It is noteworthy that the conventional data-driven MFTD with $h=0.1$ yielded some samples outperforming those from the enhanced data-driven MFTD, suggesting that the latter could further enhance its capability to discover more promising solutions through a more global search across material distributions and HF parameters.

\begin{figure}[hbt!]
\centering
\includegraphics[width=0.7\textwidth]{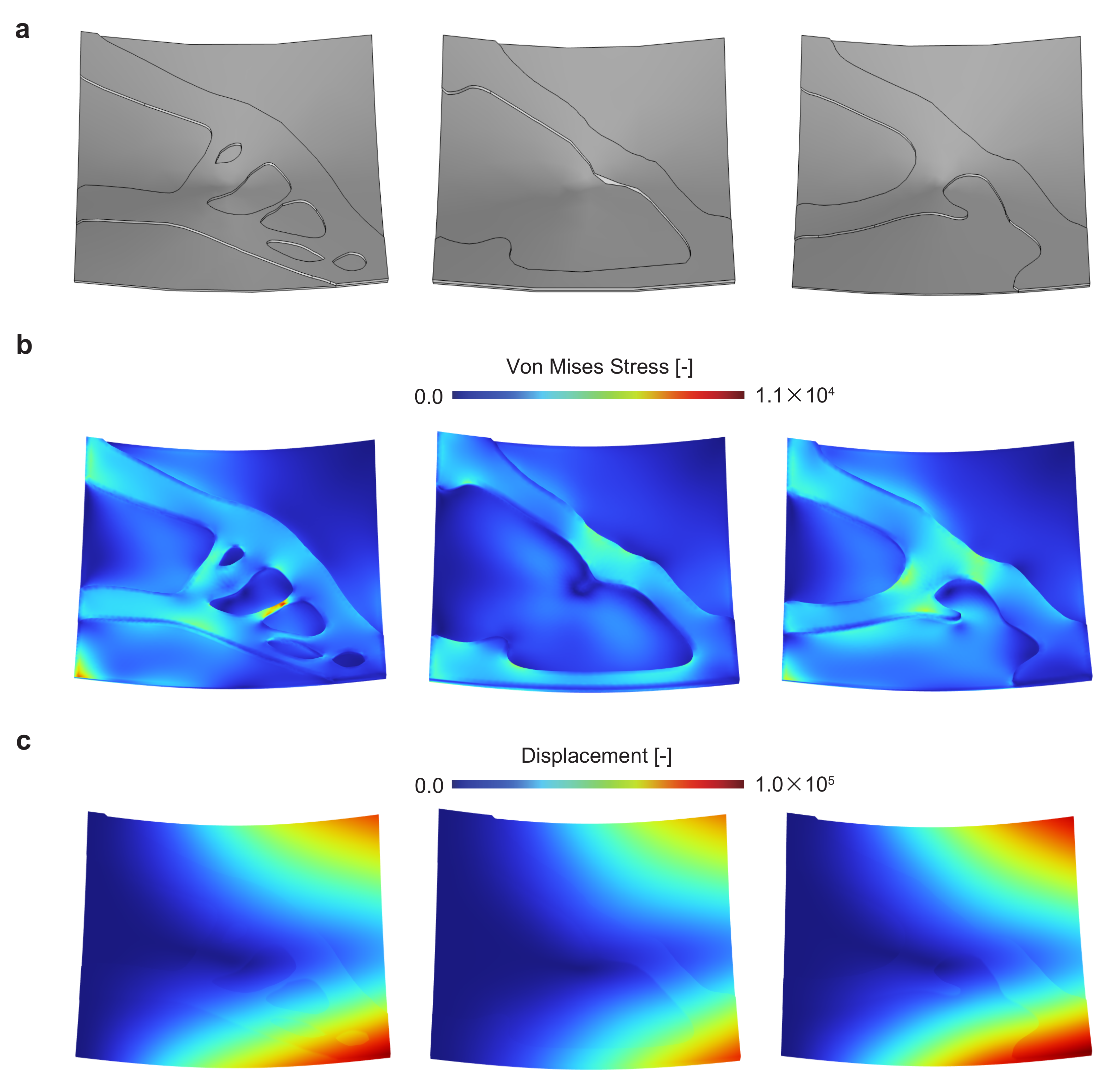}
\caption{Comparison of structural performance between stress-based initial (left), stiffness-based optimized (middle), and stress-based optimized (right) model: a) HF model; b) stress field; c) displacement field. The volume of these models are 0.015, whereas the maximum stress is $1.1\times 10^4$ (stress-based initial); $7.3\times 10^3$ (stiffness-based optmized); $7.2\times 10^3$ (stress-based optimized).}
\label{fig:init_comp_stress_ex}
\end{figure}

% 最適解・変位分布・応力分布の図
The structural performance of the optimized model obtained from LF maximum stress minimization and LF stiffness maximization, and the initial model from LF maximum stress minimization are juxtaposed in Fig.~\ref{fig:init_comp_stress_ex} on the right, middle and left, respectively.
Although all models have the same volume of $0.015$, as shown in Fig.~\ref{fig:init_comp_stress_ex} a), the initial model's reinforcement has a more complex topology with partially narrow ribs and many holes, causing a higher stress concentration as shown in Fig.~\ref{fig:init_comp_stress_ex} b).
The optimized model from both LF maximum stress minimization and LF stiffness maximization exhibit simpler topologies with uniformly thick ribs and fewer holes, leading to lower stress concentrations as shown in Fig.~\ref{fig:init_comp_stress_ex} b).
Interestingly, altough the maximum stress is lower in the order of stress-based optimized, stiffness-based optimized, and stress-based initial model, the displacement is smaller in the order of stiffness-based optimized, stress-based initial, and stress-based optimized model as shown in Fig.~\ref{fig:init_comp_stress_ex} c).
This discrepancy arises from the stress concentration depending on the topology of the reinforcement more strongly than the stiffness or deformation, clearly demonstrated by comparing the optimized models from LF stiffness maximization and LF maximum stress minimization, where the latter exhibits a lower maximum stress but a higher displacement.
It is noteworthy that the stiffness-based model has the lower displacement than both stress-based initial and optimized models, indicating that the optimized model from LF stiffness maximization reduces the stress concentration by increasing the stiffness, as expected.

%===========================================================================
\subsubsection{Effect of low-fidelity optimization formulation}
\label{subsubsec:eff_lf}
%===========================================================================

In addition to the comparison with the conventional data-driven MFTD, the proposed data-driven MFTD was further investigated by comparing the the final optimized results from two different couplings of LF optimization and HF evaluation: LF stiffness maximization and HF stress minimization, and LF maximum stress minimization and HF stress evaluation.
Initially, maximum stress minimization was approached through LF stiffness maximization coupled with HF stress evaluation.
This strategy is anticipated to indirectly yield promising solutions characterized by low maximum stress and volume by identifying samples with reduced strain energy, indicative of higher stiffness.
In a linearly elastic system, a structure with greater stiffness is expected to undergo less deformation under identical loading conditions, potentially leading to decreased stress levels within the structure, given constant applied forces.
While LF optimization for stiffness maximization is an indirect method, it offers the benefit of generally lacking non-linearity and allows for the derivation of material distributions without the need for relaxation or approximation methods.
In the second coupling, maximum stress minimization was tackled using LF approximated maximum stress minimization and HF stress evaluation, a more direct approach.
This method is expected to achieve promising solutions less indirectly through the objective functional of stress approximated by a $p$-norm, which, while still indirect, is more straightforward than stiffness maximization, where the objective functional is strain energy.
A more direct LF optimization could potentially enable the samples that inherit the LF results to attain higher performance after iterative processes.

\begin{figure}[hbt!]
\centering
\includegraphics[width=.9\textwidth]{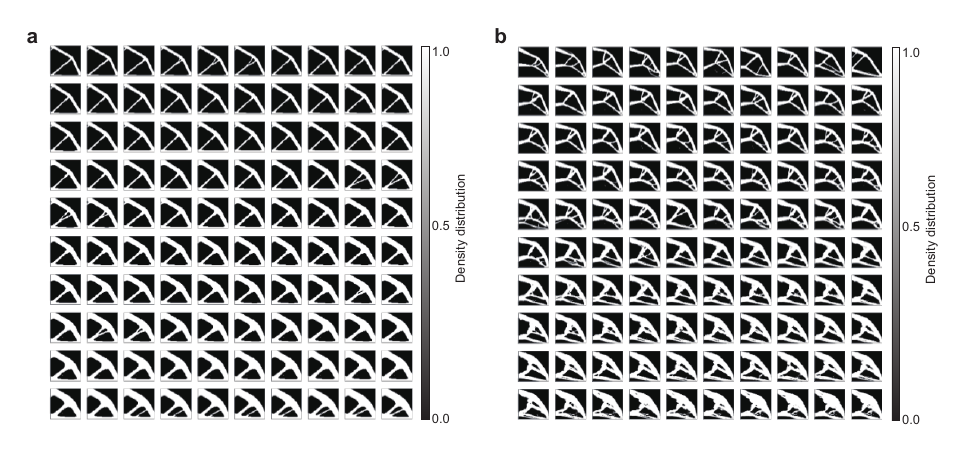}
\caption{Comparison of material distribution with $ V_{\text{max}} \in [0.2,0.5]$ from a) LF stiffness maximization and b) LF maximum stress minimization.}
\label{fig:LF_results_stress}
\end{figure}

% LFの結果: LF stiffness
The material distributions resulting from the LF stiffness maximization and LF maximum stiffness minimization are depicted in Fig.~\ref{fig:LF_results_stress}, organized according to the maximum volume constraint $V_{\text{max}}$ ranging from $0.2$ to $0.5$.
The topology of these outcomes with LF stiffness maximization closely resembles those derived from the stiffness maximization, as illustrated in Fig.\ref{fig:LF_results_stiffness}.
This resemblance is consistent with the slight distinction between two LF stiffness maximization problems, where the loading area is marginally broader for the maximum stress minimization problem.

% LFの結果: LF stress
Fig.~\ref{fig:LF_results_stress} b) reveals a wide variety of topologies can be observed in the scenario with LF maximum stress minimization, even with results based on similar volume constraints showing significant differences.
This diversity stems from the strong non-linearity inherent in the approximated maximum stress minimization problem, leading to a tendency for the optimization process to become ensnared in local optima.
Furthermore, these local optima often feature small holes and pronounced steep curves, contributing to stress concentration.
Such characteristics, arising from the approximated maximum stress minimization, present substantial opportunities to reduce the true objective functional, specifically, the original von Mises stress.

\begin{figure}[hbt!]
\centering
\includegraphics[width=0.5\textwidth]{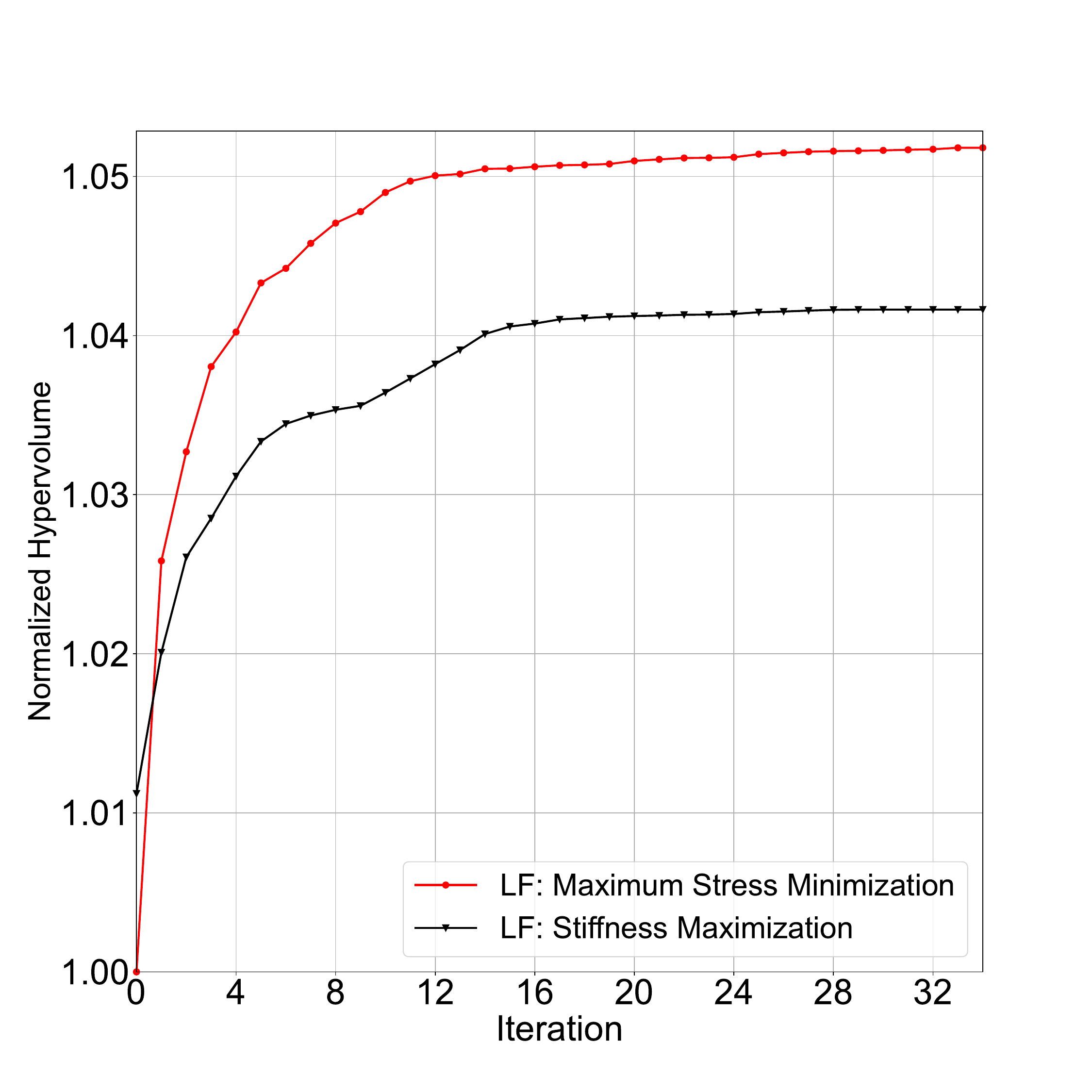}
\caption{Comparison of iteration history of hypervolume indicator between LF stiffness maximization and LF maximum stress minimization.}
\label{fig:hypervolume_history_stress}
\end{figure}

The evolution of the hypervolume indicator for the scenarios using LF stiffness maximization and LF maximum stress minimization is depicted in Fig.~\ref{fig:hypervolume_history_stress}, following normalization against the initial hypervolume indicator from the LF maximum stress minimization scenario.
The hypervolume indicators were computed using $r_{\text{hv}} = [J_1, J_2] = [14000, 0.06]$, with $J_1$ and $J_2$ representing maximum stress and volume, respectively.
The optimization concluded at iteration $i=30$ for LF stiffness maximization and at iteration $i=34$ for LF maximum stress minimization, upon the relative error of the hypervolume indicator falling below the threshold $\varepsilon_{\text{HV}}=1.0\times10^{-5}$.
The improvement in the hypervolume indicator was noted as $4.2 \%$ for LF stiffness maximization and $5.2 \%$ for LF maximum stress minimization, relative to the initial value from LF maximum stress minimization.
It is noteworthy that the initial hypervolume for the LF stiffness maximization scenario was higher than that for the LF maximum stress minimization.
This discrepancy arises because LF maximum stress minimization yielded topologies that were qualitatively sound but included numerous stress-concentrating features due to its pronounced non-linearity.
Once the iterative evolutionary process commenced, the hypervolume indicator for the LF maximum stress minimization scenario immediately exceeded that of the LF stiffness maximization scenario, improving the qualitatively sound topologies to push the Pareto front towards the utopia point in the objective space.
Furthermore, MC-VAE facilitated this rapid advancement by generating samples characterized by fewer holes and smoother material boundaries, effectively blending the features of all samples in the dataset.

%sssssssssssssssssssssssssssssssssssssssssssssssssssssssssssssssssssssssssssssss
\section{Conclusions}
\label{sec:con}
%sssssssssssssssssssssssssssssssssssssssssssssssssssssssssssssssssssssssssssssss
We proposed a data-driven MFTD framework that expands the search space for HF parameters, previously constrained to a single constant value across the whole optimization, aiming to reduce user dependency on finding optimal solutions.
This framework leverages a multi-channel image representation of material distributions and HF parameters, enabling the exploration of both material distributions and HF parameters in a single optimization process.
The effectiveness of this framework is demonstrated through two numerical examples: stiffness maximization and maximum stress minimization problems.

% Stiffness maximization
%The proposed framework was first applied to a stiffness maximization problem, a fundamental problem in structural design, to evaluate its efficacy in navigating both material distributions and HF parameters, i.e., thickness value for each solution in this case.
%The outcomes of this framework were then benchmarked against reference results obtained using a set of LF and HF parameters uniformly sampled from their target ranges.
The first example for stiffness maximization revealed that the proposed framework can successfully identify promising solutions with structural performance that is superior to that of the reference solutions by efficiently exploring both material distributions and HF parameters.
Additionally, within the data-driven MFTD's crossover process for multi-channel images, the implementation of an oversampling technique significantly improved the overall search performance by enriching the dataset that includes samples with underrepresented HF parameters.

% Maximum stress minimization: in-plane loads
%In the second example, the proposed framework was tested on a maximum stress minimization problem, a problem inherently difficult to solve directly due to its pronounced non-linearity.
%This test aimed to assess the framework's performance in handling topology optimization problems characterized by strong non-linearity.
%Additionally, the conventional data-driven MFTD was applied to this problem using four distinct HF parameters to compare its capability to efficiently navigate the solution space by adjusting both material distributions and HF parameters with the proposed data-driven MFTD.
%Furthermore, solutions for the maximum stress minimization were derived through a straightforward strategy via LF maximum stress minimization and HF stress evaluation, as well as another strategy via LF stiffness maximization and HF stress evaluation, with the latter being more indirect yet widely utilized in mechanical engineering practices.

Within the second example, the proposed framework demonstrated its ability to identify promising solutions for maximum stress minimization problems, even in the presence of significant non-linearity.
It showed a distinct advantage in globally searching for optima across material distribution and HF parameter spaces compared to the traditional data-driven MFTD method employing a SC-VAE.
The proposed framework significantly reduces computational demands and user reliance on finding optimal solutions by enabling a more comprehensive optimization in a single run, unlike traditional methods that require extensive parametric studies with varying HF parameters.

Additionally, solutions for the maximum stress minimization were derived through two distinct combinations of LF optimization and HF evaluation, where one is more straightforward and the other is more indirect yet widely utilized in mechanical engineering practices.
The latter led to more promising solutions by refining the initial dataset derived from a more direct LF topology optimization approach, which suggests that the synergy between LF topology optimization and HF evaluation critically influences the search efficiency for optimal solutions.

%Employing a strategy combining LF stiffness maximization with HF stress evaluation yielded suitable solutions by targeting samples with minimal deformation, thereby reducing maximum stress.
%Conversely, integrating LF maximum stress minimization with HF stress evaluation led to more promising solutions by refining the initial dataset derived from a more direct LF topology optimization approach.
%This outcome suggests that the synergy between LF topology optimization and HF evaluation critically influences the search efficiency for the optimal solution.

% Furture work & research direction
The HF parameter was set as a constant value across the design domain for simplicity in this study, which theoretically can be expanded to a distribution across the design domain, based on the MC-VAE's capability to process the input data stored as pixel values in multi-channel images.
We plan to address the concurrent optimization of thickness and material distributions in our future work, utilizing its further potential to optimize a distribution of HF parameters in addition to conventional material distributions.

\section*{Declaration of competing interest}
The authors declare that they have no known competing financial interests or personal relationships that could have appeared to influence the work reported in this paper.

\section*{Acknowledgements}
This work was supported by JSPS KAKENHI, Grant Number 21J22284 and 23H03799.

\appendix
%aaaaaaaaaaaaaaaaaaaaaaaaaaaaaaaaaaaaaaaaaaaaaaaaaaaaaaaaaaaaaaaaaaaaaaaaaaaaaaa
\section{Design Domain Mapping}
\label{ap:ddm}
%aaaaaaaaaaaaaaaaaaaaaaaaaaaaaaaaaaaaaaaaaaaaaaaaaaaaaaaaaaaaaaaaaaaaaaaaaaaaaaa

\begin{figure}[H]
    \centering
    \includegraphics[width=0.6\textwidth]{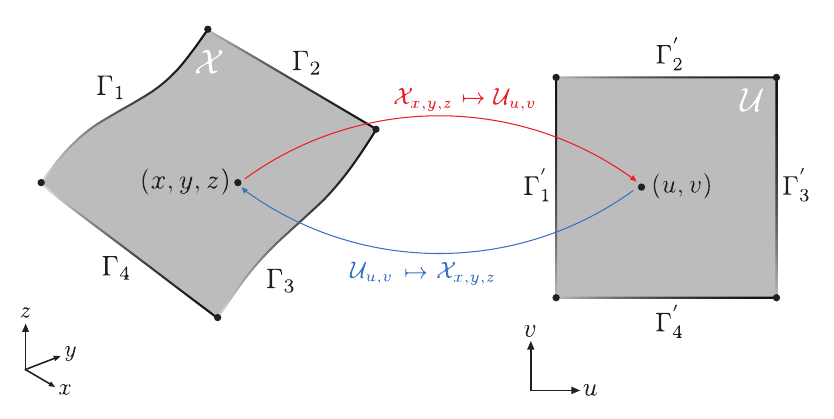}
    \caption{The concept of design domain mapping with notations.}
    \label{fig:mapping_concept}
\end{figure}

Let $\boldsymbol{x}_i = (x_i,y_i,z_i)$ represent the 3D position of the $i$-th node within a mesh patch $\mathcal{X}$, and $\boldsymbol{u}_i = (u_i,v_i)$ denote the 2D position (parameter value) of the corresponding node within the 2D unit plane mesh $\mathcal{U}$, as illustrated in Fig.~\ref{fig:mapping_concept}. The boundaries of the original mesh are divided into $\Gamma_1-\Gamma_4$, corresponding to the boundaries of the unit plane mesh denoted as $\Gamma_1^{'}-\Gamma_4^{'}$, respectively.

% Mapping notation

The boundary conditions of the unit plane mesh are divided as follows:
\begin{equation}
    \begin{aligned}
        \Gamma_1^{'} & := \{(u,v)|u=0,0\leq v\leq1\},\\
        \Gamma_2^{'} & := \{(u,v)|0\leq u\leq1,v=0\},\\
        \Gamma_3^{'} & := \{(u,v)|u=1,0\leq v\leq1\},\\
        \Gamma_4^{'} & := \{(u,v)|0\leq u\leq1,v=1\}.
    \end{aligned}
\end{equation}
\noindent
Then the mapping $\mathcal{X}_{x,y,z} \mapsto \mathcal{U}_{u,v}$ is computed by solving the following Laplace equation in the original mesh:
\begin{equation}
    \begin{aligned}
        \frac{\partial^2 u}{{\partial x}^2} + \frac{\partial^2 u}{{\partial y}^2} + \frac{\partial^2 u}{{\partial z}^2} & = 0,\\
        \frac{\partial^2 v}{{\partial x}^2} + \frac{\partial^2 v}{{\partial y}^2} + \frac{\partial^2 v}{{\partial z}^2} & = 0,
    \end{aligned}
\end{equation}
\noindent under the following Dirichlet boundary conditions:
\begin{equation}
    \begin{aligned}
        (u,v) & = (0, t_1 (x,y,z)) \quad \text{on} \quad \Gamma_1,\\
        (u,v) & = (t_2 (x,y,z), 1) \quad \text{on} \quad \Gamma_2,\\
        (u,v) & = (1, t_3 (x,y,z)) \quad \text{on} \quad \Gamma_3,\\
        (u,v) & = (t_4 (x,y,z), 0) \quad \text{on} \quad \Gamma_4.
    \end{aligned}
\end{equation}
\noindent
Herein, $t_1$ denotes the arc length ratio of $\Gamma_1^{x,y,z}$ and $\Gamma_1$, where $\Gamma_1^{x,y,z}$ denotes the boundary from the connection point of $\Gamma_4$ and $\Gamma_1$ to a point $(x,y,z)$ on $\Gamma_1$.
Solving the Laplace equation indicates the mapping minimizes the Dirichlet energy for a smooth map $\boldsymbol{u}$ from a differential surface patch $\mathcal{X}$ to its image $\mathcal{U}$, defined as follows:
\begin{equation}
    E_D = \frac{1}{2} \int_{\mathcal{X}} |\nabla \boldsymbol{u}|^2 dA,
    \label{eq:dirichlet_energy}
\end{equation}
\noindent
where $\nabla \boldsymbol{u}$ means the unit mesh $\mathcal{U}$'s gradient.
This projection distorts the mesh to accurately represent the surface's three-dimensional contours on a two-dimensional plane, ensuring all mesh elements fit within a unit square boundary.
It can be said that DDM is formulated as a \textit{harmonic mapping} \cite{eells1964harmonic}, which minimizes distortion in the mapping by mapping each node on the surface so that the Dirichlet energy is minimized given fixed boundary conditions.
Harmonic mapping can smoothly map the interior of the surface while maintaining specified boundary conditions.

Harmonic mapping is equivalent to \textit{conformal mapping} \cite{gu2012numerical} when applied under Dirichlet boundary conditions that secure all boundaries within the design domain.
The primary goal of conformal mapping is to maintain local angles throughout the mapping process by minimizing the \textit{conformal energy}, which includes the Dirichlet energy as defined in Eq.~(\ref{eq:dirichlet_energy}), alongside the area of the design domain.
The conformal mapping has been widely used in topology optimization to realize the topology optimization problem on complex surfaces \cite{vogiatzis2018computational,huo2022topology,ye2019topology}.
Given that the design domain's area remains constant under fixed boundary conditions, minimizing conformal energy effectively reduces the Dirichlet energy, thereby decreasing distortion in the mapping.
This rationale underpins the selection of the conformal mapping algorithm, favored for its straightforward implementation in Python, the primary development environment for this study.
The implementation of the conformal mapping function leverages the IGL Python library based on several works \cite{desbrun2002intrinsic,bruno2002least,mullen2008spectral}.

  \bibliographystyle{elsarticle-num} 
  \bibliography{reference}


%\begin{thebibliography}{00}


%\bibitem{}

%\end{thebibliography}
\end{document}